\documentclass[11pt]{article}

\usepackage{amsmath,amssymb}
\usepackage{graphicx}

\newtheorem{thm}{Theorem}[section]

\newtheorem{obs}{Observation}

\newtheorem{lem}[thm]{Lemma}
\newtheorem{prop}[thm]{Proposition}

\newenvironment{proof}{\noindent\bof{Proof:}}{{\hfill \rule{1ex}{1ex}}\\}
\newenvironment{oproof}[1][]{\noindent\bof{Proof of observation {#1}}\ensuremath{\blacktriangleright}}{\ensuremath{\blacktriangleleft}\\}

\newcommand{\bof}[1]{\text{\textup{\textbf{#1}}}}

\newcommand{\dimh}{\ensuremath{\textrm{\textup{{dim}}}_H}}

\newcommand{\diam}{\ensuremath{\,{\textrm{\textup{diam\,}}}}}
\newcommand{\arcdiam}{\ensuremath{\,{\textrm{\textup{arc-diam}}}}}
\newcommand{\dist}{\ensuremath{\,{\textrm{\textup{dist\,}}}}}

\newcommand{\spt}{\ensuremath{{\textrm{\textup{spt\,}}}}}
\newcommand{\gr}{\ensuremath{{\textrm{\textup{gr\,}}}}}

\newcommand{\ip}[2]{\ensuremath{\langle {#1},\, {#2} \rangle }}
\newcommand{\R}{\ensuremath{\mathbb{R}}}
\newcommand{\N}{\ensuremath{\mathbb{N}}}
\newcommand{\Q}{\ensuremath{\mathbb{Q}}}
\newcommand{\Real}{\R}
\newcommand{\Nat}{\N}

\newcommand{\dZero}{d_0}
\newcommand{\dOne}{d_1}
\newcommand{\dTwo}{\alpha_1}
\newcommand{\dThree}{c_1}
\newcommand{\dFour}{d_2}
\newcommand{\dFive}{\frac{2}{5}\frac{\alpha_0}{d_+}}
\newcommand{\dEight}{c_2}

\begin{document}
\title{The Hausdorff dimension of the visible sets of connected compact sets\footnote{2000 Mathematics Subject Classification: Primary 28A80.
Secondary 28A78, 31A15.}}
\author{Toby C O'Neil}
\date{Faculty of Mathematics and Computing, The Open University,\\ Walton Hall, Milton Keynes, MK7 6AA\\
\texttt{t.c.oneil{\makeatletter @\makeatother}open.ac.uk}\\
\textsc{Draft: }\today} \maketitle

\begin{abstract}
    For a compact set \(\Gamma\subset\R^2\) and a point \(x\), we
    define the visible part of \(\Gamma\) from \(x\) to be the set
    \[\Gamma_x=\{u\in \Gamma: [x,u]\cap \Gamma=\{u\}\}.\]
    (Here \( [x,u]\) denotes the closed  line segment joining \(x\) to
    \(u\).)

    In this paper, we use energies to show that if \(\Gamma\) is a compact connected set of
    Hausdorff dimension larger than one, then for (Lebesgue) almost every
    point \(x\in\R^2\), the Hausdorff dimension of \(\Gamma_x\) is strictly
    less than the Hausdorff dimension of \(\Gamma\). In fact, for almost every x,
    \[\dimh (\Gamma_x)\leq \frac{1}{2}+\sqrt{\dimh (\Gamma)-\frac{3}{4}}.\]

    We also give an estimate of the Hausdorff dimension
    of those points where the visible set has dimension larger
    than \(\sigma+\frac{1}{2}+\sqrt{\dimh (\Gamma)-\frac{3}{4}}\) for \(\sigma>0\).
\end{abstract}

\section{Introduction}

Given a  subset \(E\) of the plane,
Urysohn~\cite{urysohn,urysohn2} defined the notion of linear
accessibility for a point \(p\in E\): \(p\) is \emph{linearly
accessible} if there is a non-degenerate line segment \(L\) that
only meets \(E\) at the point \(p\). In a sequence of papers,
Nikodym~\cite{nikodym1,nikodym2,nikodym3} investigated the
relationship between the set theoretic complexity of \(E\) and the
set of linearly accessible points.

In this paper, we consider those points of a compact connected set
\(\Gamma\) set that are linearly accessible from a given fixed
point \(x\) and investigate the relationship between the
(Hausdorff) dimensions of the compact set and its linearly
accessible part from \(x\) for Lebsgue almost all
\(x\in\R^2\setminus \Gamma\). Denoting \(\Gamma_x\) to be the
points of \(\Gamma\) that are linearly accessible from \(x\), it
is clear that \(\dimh (\Gamma_x)\leq \dimh (\Gamma)\) for all
\(x\in\ \R^2\setminus\Gamma\).
What is perhaps surprising though
is that for most points there is a drop in dimension.

Proceeding more formally, if for a compact set in the plane,
\(K\), and \(x\in\R^2\) we  define the visible part of \(K\) from
\(x\) by
\[K_x=\{u\in K: [x,u]\cap K=\{u\}\},\]
where \( [x,u]\) denotes the closed line segment joining \(x\) to
    \(u\), then our results may be summarised as follows.

\begin{thm}\label{thmresult2}
    If \(\Gamma\subset\R^2\) is a compact connected set with \(\dimh
    (\Gamma)>1\),
    then for (Lebesgue) almost all \(x\in \R^2\),
    \[\dimh (\Gamma_x)\leq \frac{1}{2}+\sqrt{\dimh (\Gamma)-\frac{3}{4}}.\]
\end{thm}

This follows directly from the theorem that we prove in this
paper.

\begin{thm}\label{thmresult}
    Let \(\Gamma\subset\R^2\) be a compact connected set with \(\dimh
    (\Gamma)>1\).
    Then for \(\frac{1}{2}+\sqrt{\dimh (\Gamma)-\frac{3}{4}}<s\leq \dimh (\Gamma)\),
    \[\dimh \{x\in \R^2: \dimh (\Gamma_x)>s\}\leq \frac{\dimh (\Gamma)-s}{s-1}.\]
\end{thm}

In an earlier paper~\cite{jjmo}, it was shown that for a
particular class of compact connected sets (namely quasicircles),
whenever \(x\) lies outside the set, \(\dimh (\Gamma_x)=1\). Since
quasicircles can have dimension arbitrarily close to 2, and for
connected sets of positive dimension, \(\dimh (\Gamma_x)\geq 1\)
whenever \(x\not\in \Gamma\), it follows that, unless the optimal
upper bound for \(\dimh (\Gamma_x)\) is one, there is no general
result concerning the lower bound of \(\dimh (\Gamma_x)\) beyond
the trivial estimate.

There are many possible directions for future work. Despite the
fact that the upper bound given in Theorem~\ref{thmresult2} is the
golden-ratio for \(\dimh (\Gamma)=2\), there is no good reason to
believe that this bound is optimal, since the proof we give in
this paper uses at least one sub-optimal estimate. It would be
interesting to know the correct upper bound. Our method of proving
Theorem~\ref{thmresult} relies in an essential way on the
properties of connected sets in the plane, and it is unclear
whether a similar result could hold in higher dimensions. Whether
a dimension drop will occur for totally disconnected sets is also
unclear: in~\cite{jjmo}, it is shown that, for the cross-product
of a Cantor set with itself in the plane, there is a dimension
drop (to 1), provided that the original Cantor set has Hausdorff
dimension sufficiently close to 1.

I would like to thank Paul MacManus, Pertti Mattila and David
Preiss for useful discussions during the writing of this paper,
and Marianna Cs\"ornyei for her useful comments on a preliminary
draft of the paper.

\section{Background results and preliminary estimates}
In this section we summarise the main definitions and results that
we use.

Most of the time we shall be working in the plane, \(\R^2\),
endowed with the usual norm, \(|\cdot|\) and inner product
\(\ip{\cdot\,}{\cdot}\). We let \(e_1\) and \(e_2\) denote the
usual basis vectors in \(\R^2\) and set \(x^\wedge =x/|x|\) for
\(x\not=0\), and \(x^\perp=\ip{x}{e_1}e_2-\ip{x}{e_2}e_1\) for
\(x\in\R^2\).  For \(x\in \R^2\) and \(A\subseteq\R^2\), define
\(\arcdiam_x (A)\) to be the angle (in radians) subtended by the
smallest arc in the circle \(\{u:|x-u|=1\}\) that contains the
radial projection of \(A\) onto this circle. (If \(x\in A\), then
\(\arcdiam_x (A)=2\pi\).)

For a subset \(A\) of the plane and \(r>0\), let
\[B(A,r)=\{y\in\R^2: \mbox{There is \(x\in A\) with }|y-x|\leq r\}\]
and, in a slight abuse of notation, let \(B(x,r)=B(\{x\},r)\), the
usual closed ball of centre \(x\) and radius \(r\).

Let \(X\) be a Polish space. (That is, \(X\) is a complete,
separable, metrisable topological space.) A sub-additive,
non-negative set function \(\mu\) on \(X\) is a Radon measure if
it is a Borel measure (all Borel sets are \(\mu\)-measurable) for
which all compact sets have finite measure and both
\[\mu (U)=\sup\{\mu (K): K\subset U\mbox{, \(K\) is compact}\}\mbox{, for open sets }U\]
and
\[\mu (A)=\inf\{\mu (U): A\subset U\mbox{, \(U\) is open}\}\mbox{, for }A\subseteq X.\]
We denote the set of Radon measures on \(X\) by \(\mathcal{M}
(X)\).

We let \(\sigma(\bof{A}(X))\) denote the \(\sigma\)-algebra
generated by the analytic subsets of  \(X\), we suppress mention
of \(X\) when this is clear from the context. If \(\mu\) is a
Radon measure on this space then all sets in \(\sigma (\bof{A})\)
are \(\mu\)-measurable. See~\cite[(21.10)]{kechris}.

For \(s\in\R\) and \(A\subseteq X\), we define
\begin{align*}
    \lefteqn{\mathcal{M}^s (A)}\\
    &=\{\nu\in\mathcal{M}(X) : \nu (A)>0 \mbox{ and } \nu (B(x,r))\leq r^s \mbox{ for }
    x\in X,\, 0<r\leq 1\}.
\end{align*}
 If \(\mu\) is a Radon measure on the plane and
\(s\in\R\), then \(I_s (\mu)\) denotes the \(s\)-energy of \(\mu
\) given by
    \[I_s (\mu)=\iint |x-y|^{-s}\, d\mu (x)\,d\mu (y) .\]

The Hausdorff dimension of a set is defined in the usual way via
Hausdorff measures, see~\cite{falconer1,federer,mat1,rogers}. The
following theorem summarises some useful equivalent ways of
finding the Hausdorff dimension of a set.

\begin{thm}\label{thmhdim}
    Let \(A\) be an analytic subset
    of a Euclidean space, \(\R^n\). Then
    \begin{align*}
        \dimh(A) &=\sup\{s\in\R :\mathcal{M}^s (A)\not=\emptyset\}\\
        &= \sup\{s\in\R :\mbox{There is }\mu\in\mathcal{M}(\R^n)\mbox{ with }\nu (A)>0\mbox{ and }
         I_s (\mu) <\infty\}\\
         &=\sup\{\dimh (K): K\subseteq A\mbox{ and }K\mbox{ is
         compact}\}.
    \end{align*}
\end{thm}
\begin{proof}
    See~\cite[Theorem~6.4]{falconer1} together with \cite[2.10.48]{federer} or \cite[Theorem~57]{rogers}.
\end{proof}

%\subsection{Elementary geometric estimates}

We record some simple geometric estimates for future use. For
\(x\in\R^2\), \(d_-,d_+\in\R^+\), let
\(A(x,d_-,d_+)=B(x,d_+)\setminus B(x,d_-)\), a half-open annulus.

\begin{lem}
    \label{lem-angle}
    Let \( 0<d_- \leq d_+\) with \(d_-\leq 1\) and let
    \(a\in\R^2\setminus\{0\}\) and \(E\subseteq A(0,d_-,d_+)\) be compact.
    Suppose that
    \(|a|\leq \tfrac{1}{2}d_-\) and let
    \(\alpha=\min\{|\ip{p}{a^\perp}/\ip{p}{a}|:p\in E\}\). If \(\alpha\leq
    1\), then for all \(p\in E\)
    \[  \frac{1}{2}\leq \frac{\ip{p-a}{p+a}}{|p-a||p+a|}\leq 1- \frac{9}{17d_+^2}(|a|\alpha)^2.\]
\end{lem}

\begin{proof}
    For \(p\in E\),
    \[\ip{p-a}{p+a}=|p|^2-|a|^2\]
    and
    \[|p-a|^2|p+a|^2=(|p|^2+|a|^2)^2-4\ip{p}{a}^2.\]
    If \(A=\ip{p}{a^\perp}/\ip{p}{a}\), then
    \(1+A^2 =\frac{|p|^2|a|^2}{\ip{p}{a}^2}\), and so
    \[|p-a|^2|p+a|^2=\frac{(|p|^2-|a|^2)^2}{1+A^2}\left(1+
    \left(\frac{|p|^2+|a|^2}{|p|^2-|a|^2}\right)^2 A^2\right).\]
    (If \(|A|=+\infty\), then read the formula as
    \(|p-a|^2|p+a|^2=(|p|^2+|a|^2)^2\).)
    Hence
    \[\frac{\ip{p-a}{p+a}}{|p-a| |p+a|}
    =\sqrt{\frac{1+A^2}{1+(1+\mu)^2A^2}}=\sqrt{1-\frac{\mu (2+\mu)
    A^2}{1+(1+\mu)^2 A^2}},\tag{*}\]
    where
    \[2\frac{|a|^2}{d_+^2}\leq2\left(\frac{|a|}{|p|}\right)^2\leq \mu=2\frac{(|a|/|p|)^2}{1-(|a|/|p|)^2}\leq \frac{8}{3}\left(\frac{|a|}{|p|}\right)^2\leq \frac{2}{3}.\]

    It is easy to see that for \(p\in E\),~(*)
    is maximised when \(A=|\ip{p}{a^\perp}/\ip{p}{a}|=\alpha\).

    However
    \[(1-x)^{\frac{1}{2}}\leq 1-\tfrac{1}{2}x\mbox{, for \(0\leq x\leq 1\)},\]
    and so, since \(\frac{\mu (2+\mu)\alpha^2}{1+(1+\mu)^2
    \alpha^2}= 1-\frac{1+\alpha^2}{1+(1+\mu)^2\alpha^2}\leq 1\),
    and since
    \(\mu\leq \frac{2}{3}\),
    \[\frac{\ip{p-a}{p+a}}{|p-a| |p+a|}\leq 1-\frac{1}{2}\left(\frac{\mu (2+\mu)\alpha^2}{1+(1+\mu)^2
    \alpha^2}\right)\leq 1-\tfrac{9}{34}\mu \alpha^2\leq 1-\tfrac{9}{17}
    \left(\frac{|a|\alpha}{d_+}\right)^2.\]
    The lower bound follows from recognising that~(*) is minimised when
    \(p=d_- a^\perp /|a|\).
\end{proof}

For \(x\in\R^2\), \(u\in\R^2\setminus\{0\}\) and \(\sigma> 0\),
let
\[V(x,u,\sigma)=\{y\in\R^2
:|\ip{y-x}{u^{\perp}}|<\sigma\ip{y-x}{u}\},\] the open cone with
vertex \(x\), direction \(u\) and opening \(\sigma\). The next
lemma gives a lower bound on the distance of  a point in a
particular subregion of a cone from the vertex.

\begin{lem}\label{lem-intercone}
    Let \(p\in\R^2\setminus\{0\}\) and \(\sigma,\tau >0\). If
    \[u\in V(0,p,\sigma)\setminus V(p,-p,\tau),\]
    then
    \[\ip{u-p}{p^{\wedge}}\geq-\frac{\sigma}{\sigma+\tau}|p|.\]
\end{lem}

\begin{proof}
    Suppose that \(u\in V(0,p,\sigma)\setminus V(p,-p,\tau)\),
    then
    \[\ip{u-p}{p^\wedge}\geq \ip{q-p}{p^\wedge}\]
    where
    \[q=\mu (p+\sigma p^\perp )=p+\lambda (-p+\tau p^\perp ),\]
    for some \(\mu,\lambda >0\).
    Calculating \(\ip{q}{p^\perp}\) gives
    \[\mu=\lambda\frac{\tau}{\sigma}\]
    and substituting for \(\mu\) in \(\ip{q}{p}\) gives
    \[\lambda=\frac{\sigma}{\sigma+\tau}.\]
    Hence
    \[\ip{q-p}{p^\wedge}\geq -\frac{\sigma}{\sigma+\tau}|p|,\]
    as required.
\end{proof}

\subsection{Elementary measure estimates}

We now prove some estimates concerning the geometric distribution
of mass for Radon measures in the plane.

 We start by recording a
simple mass estimate.
\begin{lem}\label{lem-massnuxfx}
    Fix \(s>0\) and  \(0<d_-\leq \frac{1}{2}d_+ \). Let \(\nu\) be a Radon measure such that for all
    \(u\in\R^2\) and \(r>0\), \(\nu(B(u,r))\leq r^s\).
    Suppose that \(x\in \R^2\) and \(V\subseteq \R^2\), then
    \[\nu (V\cap A(x,d_-,d_+))\leq c\arcdiam_x (V\cap A(x,d_-,d_+))^{s-1},\]
    for some fixed positive constant \(c\)
    depending only on \(d_-,d_+\) and \(s\).
\end{lem}
\begin{proof}
    We may suppose that \(x=0\).
    Let \(\theta=\arcdiam_0 (V)\). If \(\theta\leq 1/2\), then
    \(\theta d_+\leq d_+-d_-\) and so \(V\cap A(0,d_-,d_+)\)
    may be covered by \(1+\frac{d_+-d_-}{ d_+\theta}\)
    boxes of side \(d_+\theta\).
    Hence a simple estimate of mass gives
    \[\nu (V\cap A(0,d_-,d_+))\leq 2^{\frac{1}{2}s}(d_+\theta +d_+ -d_- )(d_+\theta)^{s-1}
    \leq 3d_+^{s} 2^{\frac{1}{2}s-1}\theta^{s-1},\]
    and the lemma follows for \(\theta\leq 1/2\). If \(\theta\geq
    1/2\), then we use the estimate that \(\nu (A(0,d_-,d_+) )\leq
    d_+^s\).
\end{proof}

We now prove a lemma on the distribution of mass for an arbitrary
measure in semi-infinite tubes. To do this we define for
\(x\in\R^2\) and \(r>0\),
\[T^+ (x,r)= \{z\in \R^2 : |p_1 (z)-p_1(x)|<r \mbox{ and }p_2 (z)>p_2 (x)\}\]
and
\[T^- (x,r)= \{z\in \R^2 : |p_1 (z)-p_1(x)|<r \mbox{ and }p_2 (z)<p_2 (x)\},\]
where \(p_1\) and \(p_2\) denote orthogonal projection onto the
\(x\)-\ and \(y\)-axis, respectively. Thus \(T^+ (x,r)\) is an
open vertical tube of width \(2r\) extending upwards from \(x\)
and \(T^- (x,r)\) is an open vertical tube of width \(2r\)
extending downwards from \(x\).

\begin{prop}\label{prop-lines}
%    Let $s>1$ and $c>0$.
    Suppose  $\nu$ is a compactly supported
    Radon measure in the plane.
%    for all $x$
%    and $0<r\leq 1$,
%    \[ \nu (B(x,r)) \leq cr^s .\]
    Then for $\xi>0$
    and $\nu $-a.e.\ $x$
    \[\liminf_{r\rightarrow 0} \frac{\nu( T^{+}(x,r))}{r^{1+\xi}}=
    \liminf_{r\rightarrow 0}\frac{\nu( T^{-}(x,r))}{r^{1+\xi}}=+\infty.\]
\end{prop}
\begin{proof}
    We give the proof for \(T^+\); the proof for \(T^-\) is
    similar.
    Without loss of generality we assume that \(\spt\nu\) lies in the
    unit square \([0,1]\times [0,1]\) and let
    \[E_\infty=\left\{x:\liminf_{r\rightarrow 0} \frac{\nu( T^{+}(x,r))}{r^{1+\xi}}=
        +\infty \right\}.\]
%    We omit the elementary verification that \(E_\infty\) is a
%    Borel set.
    Since \(\nu (\R^2\setminus\spt\nu)=0\), it is
    enough to show that \(\nu (\spt\nu\setminus E_\infty )=0\).

    For \(M\) and \(j\in\N\cup\{0\}\), let
    \[E^{M,j} =\{x\in \spt\nu : \nu ( T^+ (x,r))<M r^{1+\xi}\mbox{ for some }
    0<r\leq 2^{-j}\}.\]
    Then
    \[\spt\nu=E_\infty \cup\bigcup_{M\in\N}\bigcap_{j\in\N\cup\{0\}}E^{M,j},\]
    and
    \[E^{M,j}\subset
    \bigcup_{k\geq j}\underbrace{\{x\in\spt\nu :\nu ( T^+ (x,2^{-k}))
    <2^{1+\xi}M 2^{-k(1+\xi)} \}}_{E^{M,j,k}, \mbox{ say}}.\]
    We now estimate the \(\nu\) measure of \(E^{M,j,k}\) for \(k\geq j\in\Nat\).
    Choose \(F\subseteq E^{M,j,k}\) compact such that
    \[\nu (F)\geq \nu (E^{M,j,k})/2 .\]
%    (We omit the verification that \(E^{M,j,k}\) is a
%    Borel set for each \(M\), \(j\) and \(k\).)
    We consider the \(2^{k+2}\) columns \(C_i
    =[i 2^{-(k+2)},(i+1)2^{-(k+2)}]\times\Real\), \(i=0,\ldots,
    2^{k+2}-1\).
    For each \(i\) with \(C_i\cap F\not=\emptyset\), we choose \(x_i\in
    C_i\cap F\) to have minimum possible height above the
    \(x\)-axis, ie
    \[\dist (C_i\cap F,\R\times \{0\})=\dist (x_i,\R\times
        \{0\}).\]
    For such an \(i\),
    \[\nu (F\cap T^+ (x_i ,2^{-k})\leq \nu( T^+ (x_i ,2^{-k})) <
    2^{1+\xi}M 2^{-k(1+\xi)}.\]
    Clearly
    \[F\subseteq\bigcup_{i: F\cap C_i\not=\emptyset}F\cap C_i
    \subseteq \bigcup_{i: F\cap C_i\not=\emptyset}F\cap T^+ (x_i ,2^{-k}).\]
    And so
    \begin{eqnarray*}
        \nu (F) & \leq & \sum_{i: F\cap C_i\not=\emptyset}\nu
        (F\cap T^+ (x_i ,2^{-k})) \\
        & < & 2^{k+2}\times 2^{1+\xi}M 2^{-k(1+\xi)}\\
        & = & 2^{3+\xi}M 2^{-k\xi}.
    \end{eqnarray*}
    Hence
    \[\nu (E^{M,j,k}) < 2^{4+\xi}M 2^{-k\xi}\]
    and so
    \[\nu (E^{M,j})\leq \sum_{k=j}^\infty \nu (E^{M,j,k}) <
    \frac{2^{4+\xi}M}{1-2^{-\xi}}2^{-j\xi}.\]
    Thus,
    \[\nu\left(\bigcup_{M\in\N}\bigcap_{j\in\N} E^{M,j}\right)=0\]
    and the lemma follows.
\end{proof}

\begin{figure}[!htb]
  \centering
  \includegraphics[height=1.5in, clip=true, keepaspectratio=true]{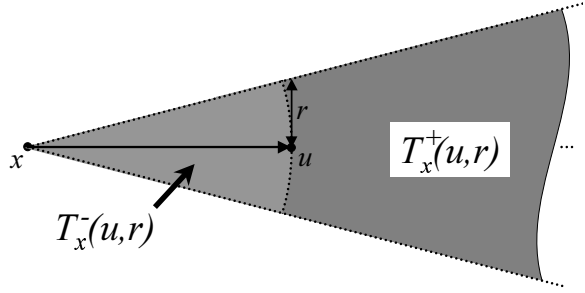}
  \caption{The radial tubes \(T_x^- (u,r)\) and \(T_x^+ (u,r)\).}\label{fig-radtubes}
\end{figure}

 For
\(x\not=u\in\R^2\) and \(r>0\), define radial tubes \(T_x^+(u,r)\)
and
 \(T_x^-(u,r)\) by
\[T_x^+(u,r)= V(x,u-x,r/d(x,u))\cap \{z\in\R^2: d(x,z)> d(x,u)\}\]
 and
\[T_x^-(u,r)=V(x,u-x,r/d(x,u))\cap \{z\in\R^2: d(x,z)<d(x,u)\},\]
 see
Figure~\ref{fig-radtubes}.

It is easy to use a bi-Lipschitz transformation to transform our
lemma about parallel tubes to one about radial tubes.

\begin{lem}
    \label{lem-radtubes}\label{lem-cones}
    Let  \(\nu\) be a compactly supported Radon measure in the plane
    and \(x\not\in\spt\nu\).
    Then for \(\xi >0\) and for \(\nu\)-a.e.\ \(u\)
    \[\liminf_{r\to 0}\frac{\nu (T^+_x (u,r))}{r^{1+\xi}}=
    \liminf_{r\to 0}\frac{\nu (T^-_x (u,r))}{r^{1+\xi}}=+\infty .\]
\end{lem}
\begin{proof}
    Since \(x\not\in\spt\nu\), there is \(\rho>0\)
    with \(B(x,\rho)\cap\spt\nu=\emptyset\).
%    Similarly, we may suppose that  \(E\) is disjoint from \(B(x,\rho)\).
    Since \(\spt\nu\) is compact, we can find some \(R>\rho\) for which
    \(\spt\nu\subset B(x,R)\).
    Moreover, by restricting and translating \(\nu\) suitably, we
    may suppose that \(\spt\nu\) is a subset of a quadrant of
    the plane with corner at \(x\), \(Q(x)\), say, intersected with the annulus \(A (x,\rho/2,R)\). It is now straightforward to find a
     transformation (namely, \(re^{i\theta}\mapsto (r,\theta )\))
    which transforms radial lines segments through \(x\) and intersecting
    this region to half-lines parallel to the \(y\)-axis. This transformation is bi-Lipschitz
    when restricted to \( Q(x)\cap A(x,\rho/2,R)\).  This gives us the situation
    described in Proposition~\ref{prop-lines} and the claim follows.
\end{proof}

This lemma allows us to show that measures with dimension larger
than one have mass far from the origin of these radial tubes for
typical points:

\begin{lem}\label{lem-Tmass}
    Let \(s>1\), \(0<r_1\leq  r_0\leq 1\)
    and \(\xi, M, d_-, c>0\), and
    \(x\in\mathbb{R}^2\). Suppose that \(\nu\) is a compactly supported Radon
    measure on the plane and \(F\subseteq E\) are compact sets in the
    plane satisfying:
    \begin{enumerate}
        \item for all \(u\in E\), \(|u-x|\geq d_-\);
        \item for all \(u\in E\) and \(0<r\leq r_0\),
        \[\nu B(u,r)\leq cr^s ;\]
        \item for all \(u\in F\) and \(0<r\leq r_1\),
        \[\nu (E\cap T_x^{\pm} (u,r))> Mr^{1+\xi}.\]
    \end{enumerate}
    Then there are constants \(r_2\in (0,r_1/\sqrt{2}]\) and \(\dZero>0\) such that for \(u\in F\) and
    \[0<r \leq r_2,\]
    \[\nu (E\cap T_x^\pm (u,r)\cap(\mathbb{R}^2\setminus A(x,|u-x|-\dZero r^{2+\xi-s},
        |u-x|+\dZero r^{2+\xi-s})))>0.\]

\end{lem}
\begin{proof}
    Let
    \[\dZero= \tfrac{M}{12}2^{-s/2}\mbox{ and }
    r_2=\min\{r_1/\sqrt{2},\tfrac{1}{2}\sqrt{3}d_-,(d_-/\dZero)^{\frac{1}{2+\xi-s}},\dZero^{\frac{1}{s-1-\xi}}\}.\]

    We give the proof for \(T_x^+ (u,r)\); the proof for \(T_x^-
    (u,r)\) is similar. By rotating and translating, we may assume that \(x=0\) and the line segment \([x,u]\) is
    on the positive \(x\)-axis. Let \(\Delta=|u-x|\geq d_-\).
\begin{figure}
    \begin{center}
        \includegraphics[keepaspectratio=true, clip=true, width=20em]{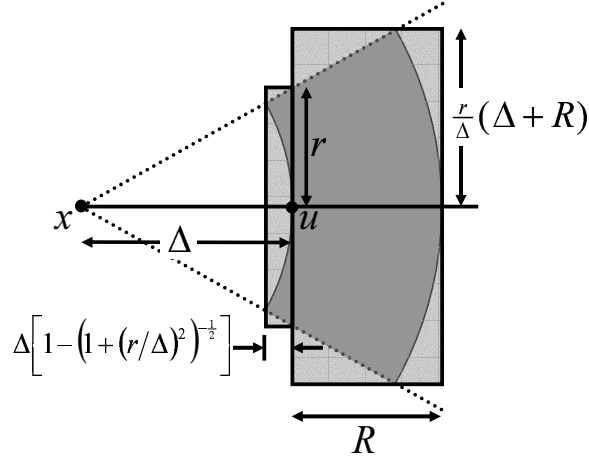}
    \end{center}
    \caption{Estimating the mass of $T^+_x (u,r)$.}\label{fig-masstubes}
\end{figure}

    Elementary geometry shows, since \(r\leq \tfrac{1}{2}\sqrt{3}d_-\leq
    \tfrac{1}{2}\sqrt{3}\Delta\) and so \(
    (1+(r/\Delta)^2)^{-\frac{1}{2}}\geq
    1-\frac{1}{2}(r/\Delta)^2\),  that
    \begin{align*}
        \lefteqn{T_x^+ (u,r)\cap A(x,\Delta-R,\Delta+R)\subseteq T_x^+ (u,r)
        \cap B(x,\Delta+R)}\\
        &\subseteq \left([\Delta-\tfrac{1}{2}r^2/\Delta
        ,\Delta]\times [-r,r]\right)\cup \left([\Delta,\Delta+R]\times [-r (1+R/\Delta ),
        r (1+R/\Delta)]\right),
    \end{align*}
    for any \(R\geq 0\). We choose \(R=\dZero r^{2+\xi -s}\).

    We estimate that \(T_x^+ (u,r)\cap (A,\Delta-R,\Delta+R)\) can
    be covered by
    \[2+(1+2R/r)(r(1+R/\Delta)+1)\]
    closed squares of side \(r\), since
    \(\tfrac{1}{2}r^2/\Delta \leq \tfrac{1}{2}\tfrac{r_2}{d_-}r<r\).
    We find that
    \[2+(1+2R/r)(r(1+R/\Delta)+1)\leq 2R/r +(3R/r)(2r+1)\leq 11 R/r , \]
    since \(r\leq R\leq \Delta\) and \(r\leq 1\).
    Hence we require at most \(11R/r\)
    balls of radius \(\sqrt{2}r\) to cover
    \(E\cap T_x^+ (u,r)\cap (A,\Delta-R,\Delta+R)\).

    So, since \(\sqrt{2}r\leq r_1\leq r_0\), we estimate that
    \begin{align*}
        \nu (E\cap T_x^+ (u,r)\cap (A,\Delta-R,\Delta+R)) &\leq (11R/r)\times
    2^{s/2}r^s\\
    &= 11\cdot 2^{\frac{s}{2}} \dZero r^{1+\xi}\\
    &<Mr^{1+\xi},
    \end{align*}
    proving the lemma.
\end{proof}

\subsection{A `two measures' estimate}

In this subsection, we investigate the interaction of two
 measures of large dimension when they are
supported on different visible sets of \(\Gamma\). The result that
we prove in this section is the crux of our method. It shows that
if two measures of large dimension are supported in different
visible sets, then they will be `disjoint' in the sense that balls
containing points from both visible sets will have small mass for
both measures. The remainder of the paper consists mainly of
trying to place ourselves in a position to use this observation.

 In
the following proposition, \(\mathcal{T}(x,y,p)\) denotes the
closed triangle with vertices \(x\), \(y\) and \(p\), and
\(H(x,y;u)\) denotes the closed  upper-half plane that has the
line segment \([x,y]\) in its boundary and \(u\) lying in its
interior.

\begin{prop}\label{prop-geom1}
    Let \(\Gamma\) be a non-empty compact connected subset of
    \(\R^2\).
    Suppose that \(s>1\), \(0<\xi<s-1\), \(0<r_1\leq r_0\leq 1\), \(0<d_-\leq d_+\) with \(d_-\leq 1\) and
    \(M>0\) are given.
    Let \(x,y\in \R^2\setminus\Gamma\) satisfy
    \[0<2|x-y|<d_-\leq \min\{d(x,\Gamma),\, d(y,\Gamma)\}
        \leq \max\{ d(x,\Gamma),\, d(y,\Gamma)\}+|\Gamma|\leq d_+ .\]
    Let \(\nu_x\) and \(\nu_y\) be Radon measures supported in \(\Gamma_x\)
    and \(\Gamma_y\) respectively and let
    \[F_x\subseteq E_x\subseteq \Gamma_x\mbox{ and } F_y\subseteq E_y\subseteq
    \Gamma_y\]
    be compact sets.
    Suppose that:
    \begin{enumerate}
        \item for all \(u\in E_x\), \(v\in E_y\) and \(0<r\leq
        r_0\) both
        \[\nu_x (B(u,r)\leq r^s\mbox{ and }\nu_y (B(v,r))\leq r^s;\]
        \item for all \(u\in F_x\), \(v\in F_y\) and \(0<r\leq
        r_1\) both
        \[\nu_x (T_x^\pm(u,r)\cap E_x)\geq Mr^{1+\xi}\mbox{ and }
        \nu_y (T_y^\pm(v,r)\cap E_y)\geq Mr^{1+\xi};\]
        \item there is \(\psi\in (0,1/2)\) such that for
        \(u\in F_x\cup F_y\),
        \[\ip{(u-x)^\wedge}{(u-y)^\wedge}\in [\tfrac{1}{2},1-\psi].\]
    \end{enumerate}
    Then there are constants \(\alpha_0, \dOne, \dThree>0\) such that for \(u\in F_x\),
    if \(0<\rho\leq \dOne\psi^{\frac{1}{2}\frac{1}{s-1-\xi}}\), then
    \begin{equation}
        \label{eqn-geom1}
        \nu_y (B(u,\rho)\cap F_y)\leq \dThree \psi^{-\frac{1}{2}\left(\frac{s-1}{2+\xi-s}
        \right)}\rho^{\frac{1+\xi}{2+\xi-s}}.
    \end{equation}

    Furthermore, if \( v\in F_y \cap B(u,\rho )\), then there is
    \[p\in [\tfrac{1}{2}(x+y),u]\cap B (u,\alpha_0 \psi^{-1/2}\rho)\]
    such that \(\mathcal{T}(x,y,p)\cap \Gamma=\emptyset\)
    and
    \[V(p,\tfrac{1}{2}(x+y)-u,\tfrac{2}{5}\psi^{\frac{1}{2}})\cap
    \Gamma\cap H(x,y;u)=\emptyset. \]
\end{prop}

Notice that the symmetry of the hypotheses in this proposition
imply that a version of~(\ref{eqn-geom1}) holds for \(u\in F_y\)
with \(\nu_y\) replaced by \(\nu_x\) and \(F_y\) replaced by
\(F_x\).

Here \(\alpha_0=60 d_+/d_-\),
\(\dOne=\min\{(r_2/\dTwo)^{2+\xi-s},\, d_-/\alpha_0\}\) and
\(\dThree=2^{5+s/2}\dTwo^{s-1}(d_+/d_-)\), where
\(\dTwo^{2+\xi-s}=(\alpha_0+1)/\dZero\) and \(\dZero,\) and
\(r_2\) are the constants determined in Lemma~\ref{lem-Tmass}.

\begin{proof}
    Suppose the conditions of the proposition are satisfied. Fix
    \[0<\rho\leq \dOne \psi^{\frac{1}{2}\frac{1}{s-1-\xi}},\]
    we must show that
    \[\nu_y (F_y\cap B(u,\rho))\leq \dThree\psi^{-\frac{1}{2}\left(\frac{s-1}{2+\xi-s}
        \right)}
    \rho^{\frac{1+\xi}{2+\xi-s}}.\]
%    Notice that this is a non-empty interval of values for \(\rho\)
%    since \(r^{s-1-\xi}\leq c_1\sqrt{\psi}\).

    If \(F_y\cap B(u,\rho )=\emptyset\), then there is nothing to
    prove. So suppose \(w\in F_y\cap B(u,\rho)\) and set
    \[e=(u-x)^\wedge,\quad f=(w-y)^\wedge\quad \mbox{and}\quad g=(u-y)^\wedge .\]
    \begin{figure}
        \begin{center}
            \includegraphics[width=0.7\textwidth, clip=true, keepaspectratio=true]{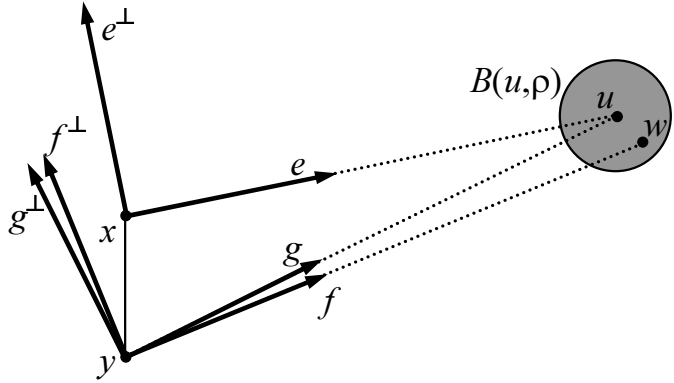}
        \end{center}
        \caption{The vectors \(e\), \(f\) and \(g\).}
    \end{figure}
    Notice that the third hypothesis of the proposition states
    \begin{equation}
        \tfrac{1}{{2}} \leq \ip{e}{g} \leq 1-\psi
        \label{eqn-cl11}
    \end{equation}
    and since
    \[\ip{e^\perp}{g}^2=1-\ip{e}{g}^2\geq 1-(1-\psi)^2=\psi (2-\psi )\geq \psi,\]
%    and \(\ip{e^\perp}{g}>0\),
    it follows that
    \begin{equation}
        |\ip{e^\perp}{g}| \geq {\psi}^{\frac{1}{2}}.\label{eqn-cl14}
    \end{equation}

    In order to prove the theorem, we make a sequence of geometric observations.
    In the first observation, we make some further estimates
    relating the angles between various of the
    vectors \(e,e^\perp, f,f^\perp,g\) and \(g^\perp\).

\begin{obs}\label{obs-one}
    If %\(|u-w|\leq \rho\) %, \(\ip{e}{g}\in [1/\sqrt{2}, 1-\psi]\) and
    \(0<\rho< \tfrac{1}{4} d_- \), then
    \begin{align}
         |\ip{f}{g^\perp}| &\leq \rho/d_- \label{eqn-cl12}\\
        \ip{f}{g}&\geq 1-2\rho/d_- \label{eqn-cl13}\\
         \ip{e}{f} &> \tfrac{1}{4}d_-/d_+ .\label{eqn-cl15}
    \end{align}
%    and if, moreover, \(\rho< \tfrac{1}{4} d_- (d_-/d_+)^2\), then \(\ip{e}{g}>0\).
\end{obs}
\begin{oproof}[\ref{obs-one}]
    For inequality~(\ref{eqn-cl12}), we use \(g^\perp =\ip{g}{f}f^\perp
    -\ip{g}{f^\perp}f\) and calculate
    \[\ip{f}{g^\perp}=0-\ip{g}{f^\perp}=-\frac{1}{|u-y|}\ip{u-y}{f^\perp}
    =-\frac{\ip{u-w}{f^\perp}}{|u-y|}.\]
    Hence \(|\ip{f}{g^\perp}|\leq \rho /d_-\).

    For inequality~(\ref{eqn-cl13}), on noting
    \[\ip{f}{g}|w-y| |u-y| =\ip{w-y}{u-y}= \ip{w-u}{u-y}+|u-y|^2,\]
    we find
    \[\ip{f}{g}=\frac{1}{|w-y|}\ip{w-u}{g}+\frac{|u-y|}{|w-y|}
    \geq -\frac{\rho}{d_-}+\left(1-\frac{|u-w|}{|w-y|}\right)\geq 1-2\rho /d_- .\]

    To verify inequality~(\ref{eqn-cl15}), note that \(w=y+|w-y|f\in B(u,\rho
    )\), and so \(w=y+(x-y)+|u-x|e+z\) for some \(|z|\leq \rho\).
    Hence
    \(|w-y|f=(x-y)+|u-x|e+z\) and
    \[|w-y|\ip{f}{e}=\ip{x-y}{e}+|u-x|+\ip{z}{e}.\]
    Now
    \[|\ip{x-y}{e}|\leq \tfrac{1}{2}d_-\leq \tfrac{1}{2}|u-x|
    \quad\mbox{and}\quad |\ip{z}{e}|\leq\rho\leq \tfrac{1}{4}d_-\leq
    \tfrac{1}{4}|u-x|.\]
    Thus
    \[|w-y|\ip{e}{f}\geq \tfrac{1}{4}|u-x|\geq\tfrac{1}{4}d_-\]
    and so \(\ip{e}{f}\geq  \tfrac{1}{4}d_-/d_+ \), as required.
%
%    Finally we note that \(|u-y|g=|w-y|f+(u-w)\), and so
%    \[\ip{g}{e}=\frac{|w-y|}{|u-y|}\ip{f}{e}+\frac{1}{|u-y|}\ip{u-w}{e}
%    \geq \tfrac{1}{4}\left(d_-/d_+\right)^2-\rho/d_- >0,\]
%    provided \(\rho <\tfrac{1}{4} d_- (d_-/d_+)^2\).
\end{oproof}

We now note that if \(z\in T_y (w,r)\), then it is also in \(T_y
(u, r') \) for \(r'\) not too much bigger than \(r\).

\begin{obs}\label{obs-two}
    If \(0< \rho \leq \tfrac{1}{4}d_- \), then % and \( |u-w|\leq \rho \)
    \[V\left(y,f,\tfrac{\rho}{d_-}\right)\subseteq V\left(y,g,4\tfrac{\rho}{d_- }\right).\]
\end{obs}

\begin{oproof}[\ref{obs-two}]
    If \(z\in V(y,f,\rho/d_-)\), then
    \begin{equation}
        |\ip{z-y}{f^\perp}|<\frac{\rho}{d_-}\ip{z-y}{f}.\label{eqn-cone1}
    \end{equation}
    Since
    \[z-y=\ip{z-y}{f}f+\ip{z-y}{f^\perp}f^\perp,\]
    we find
    \[\ip{z-y}{g^\perp}=\ip{z-y}{f}\ip{f}{g^\perp}+
    \ip{z-y}{f^\perp}\ip{f^\perp}{g^\perp}.\]
    Hence~(\ref{eqn-cl12}) implies
    \[|\ip{z-y}{g^\perp}|\leq \frac{\rho}{d_-}|\ip{z-y}{f}|+|\ip{z-y}{f^\perp}|.\]
    Thus~(\ref{eqn-cone1}) gives
    \begin{align}
        |\ip{z-y}{g^\perp}| &\leq
        \frac{\rho}{d_-}|\ip{z-y}{f}|+\frac{\rho}{d_-}\ip{z-y}{f}\nonumber\\
        &=  2(\rho/ d_-)\ip{z-y}{f}\label{eqn-cone2}
    \end{align}
    It only remains to estimate \(\ip{z-y}{f}\) in terms of
    \(\ip{z-y}{g}\). As \(f=\ip{f}{g}g+\ip{f}{g^\perp}g^\perp\),
    \[0<\ip{z-y}{f}\leq \ip{z-y}{g}\ip{f}{g}+\ip{z-y}{g^\perp}\ip{f}{g^\perp},\]
    which, on using~(\ref{eqn-cone2}) and~(\ref{eqn-cl12}), gives
    \[0<\ip{z-y}{f}\leq \ip{z-y}{g}\ip{f}{g}+\frac{2\rho}{d_-}\times
    \frac{\rho}{d_-}\ip{z-y}{f}.\]
    Rearranging and using \(0<\ip{f}{g}\leq 1\),
%    Combining this with~(\ref{eqn-cl13}),
we find
    \[\ip{z-y}{f} [1-2(\rho/d_-)^2]\leq \ip{z-y}{g}.\]
    Substituting back into~(\ref{eqn-cone2}), then gives
    \[|\ip{z-y}{g^\perp}|\leq 2(\rho /d_-)[1-2(\rho/d_-)^2]^{-1}
    \ip{z-y}{g}\]
    which, as \( \rho \leq d_- /2\), proves the claim.
\end{oproof}

\begin{obs}\label{obs-three}
    If
    \(0< \rho \leq \tfrac{1}{20}d_- %(d_-/d_+)^2
    {\psi^{1/2}}\),
    then
    \[V(x,e,\rho/d_-)\cap V(y,f,\rho /d_-)\subseteq B(u,\alpha_0
    \psi^{-1/2}\rho),\]
    where \(\alpha_0 = 60 \frac{d_+}{d_-}\).
\end{obs}

\begin{oproof}[\ref{obs-three}]
    Fix \(z\in V(y,f,\rho /d_-)\cap V(x,e,\rho /d_-)\).
    Since \(0<\rho\leq \tfrac{1}{20}d_-{\psi^{1/2}} \leq d_-/4\), observation~\ref{obs-two} implies
    \(z\in V(y,g,4\rho /d_-)\). Hence there are \(\lambda, \mu>0\)
    for which
    \[z=y+\lambda (g-bg^\perp )=x+\mu (e+ae^\perp )\]
    where \(|b|\leq 4\rho /d_-\) and \(|a|\leq \rho/d_-\). We wish
    to find an upper bound for \(|z-u|\). Now
    \[\ip{z-x}{e}=\mu\quad\mbox{and}\quad\ip{z-y}{g}=\lambda .\]
    Notice that
    \begin{equation}
        |z-u|^2=\ip{z-u}{g}^2+\ip{z-u}{g^\perp}^2=(\lambda-|y-u|)^2+b^2\lambda^2,\label{eqn-obs3-1}
    \end{equation}
    and so upper estimates for \((\lambda-|y-u|)^2\) and \(\lambda^2\) give an upper estimate for
    \(|z-u|\).

    Now
    \[\ip{z-u}{e}=\ip{y-u}{e}+\lambda(\ip{g}{e}-b\ip{g^\perp}{e})=
    \ip{x-u}{e}+\mu\]
    and so
    \[\mu=|x-u|-|u-y|\ip{g}{e}+\lambda (\ip{g}{e}-b\ip{g^\perp}{e}).\]
    Also
    \[\ip{z-u}{e^\perp}=\ip{y-u}{e^\perp}+\lambda (\ip{g}{e^\perp}-
    b\ip{g^\perp}{e^\perp})=a\mu\]
    and so
%    \[a\mu=-|u-y|\ip{g}{e^\perp}+\lambda(\ip{g}{e^\perp}-
%    b\ip{g^\perp}{e^\perp}).\]
%    Hence
    \begin{align*}
        \lefteqn{-|u-y|\ip{g}{e^\perp}+\lambda (\ip{g}{e^\perp}-
            b\ip{g^\perp}{e^\perp})}\\
        &=a|u-x|-a|u-y|\ip{g}{e}+a\lambda (\ip{g}{e}-b\ip{g^\perp}{e}).
    \end{align*}
    This rearranges to give
    \[\lambda\gamma=a|u-x|+|u-y|(\ip{g}{e^\perp}-a\ip{e}{g}),\]
    where
    \[\gamma=(1-ab)\ip{g}{e^\perp}-(a+b)\ip{e}{g}.\]

    Thus
    \begin{align*}
        \lambda-|u-y| &=\gamma^{-1}\left[a|u-x|+|u-y|(\ip{g}{e^\perp}-a\ip{e}{g}-
        \gamma)\right]\\
        &=\gamma^{-1}\left[a|u-x|+|u-y|(ab\ip{g}{e^\perp}+b\ip{e}{g})\right]\\
        &=\gamma^{-1}\left[a|u-x|+b|u-y|(a\ip{g}{e^\perp}+\ip{e}{g})\right].
    \end{align*}

    Since \(|a|\leq \rho/d_-\) and \(|b|\leq 4\rho/d_-\),
    it follows that \(|ab|\leq 1/2\) and \( |a+b|\leq 5\rho/d_-\). From
    equation~(\ref{eqn-cl14}) we  know
    \(|\ip{e^\perp}{g}|\geq\psi^{1/2}\), and so
    \[|\gamma |\geq \tfrac{1}{2}\psi^{1/2}-5(\rho/d_-)\geq \tfrac{1}{4}
    \psi^{1/2},\]
    since \(\rho\leq \tfrac{1}{20}d_-\psi^{1/2}\).

    Hence, as \(|a|\leq \rho/d_-\leq 1\),
    \[|\lambda|\leq 4\psi^{-1/2}\left[|a||u-x|+|u-y|(|\ip{g}{e^\perp}|
    +|a||\ip{e}{g}|)\right]\leq 12d_+\psi^{-1/2}\]
    and, as \(|b|\leq 4\rho/d_-\),
    \begin{align*}
        \lefteqn{|\lambda-|u-y||} \\
        &\leq 4\psi^{-1/2}\left[|a||u-x|+|b||u-y|(|a||\ip{g}{e^\perp}|
        +|\ip{e}{g}|)\right]\\
        &\leq 36 (d_+/d_-)\psi^{-1/2}\rho .
    \end{align*}
    Thus estimating \(\lambda\) in~(\ref{eqn-obs3-1}) gives
    \[|z-u|^2\leq (36 (d_+/d_-)\psi^{-1/2}\rho)^2+(48(d_+/d_-)\psi^{-1/2}\rho)^2
    ,\]
    and so
    \[|z-u|\leq 60 (d_+/d_-)\psi^{-1/2}\rho ,\]
    as required.
\end{oproof}

We now observe that there is a `large' triangle that is disjoint
from \(\Gamma\) and with a vertex close to \(u\) (and hence
\(w\)).

\begin{obs}\label{obs-four}
    If
    \[0<\rho\leq \dOne \psi^{\frac{1}{2}\left(\frac{1}{s-1-\xi}\right)},\]
    and \(r=\dTwo (\psi^{-\frac{1}{2}}\rho)^{\frac{1}{2+\xi-s}}\),
    then there is
    \[z\in V(x,e,r /d_-)\cap V(y,f,r/d_-)\]
    with
    \[\mathcal{T}(x,y,z)\cap\Gamma=\emptyset\]
    and
    \[z\in B(u,\alpha_0\psi^{-1/2}\rho).\]
\end{obs}

\begin{oproof}[\ref{obs-four}]
    We aim to find a point \(z\) which is visible from both \(x\)
    and \(y\).
    Recall that
    \[\dTwo^{2+\xi-s}=(\alpha_0+1)/\dZero\mbox{ and }
    \dOne=\min\{(r_2/\dTwo)^{2+\xi-s},\, d_-/\alpha_0\}.\]
    (The constant \(\alpha_0\)
    is given in observation~\ref{obs-three}, and \(r_2\) and \(\dZero\) are given
    in Lemma~\ref{lem-Tmass}.)

    Since \(w\in F_y\) and
    \[r=\dTwo (\psi^{-\frac{1}{2}}\rho)^{\frac{1}{2+\xi-s}}\leq \dTwo \dOne^{\frac{1}{2+\xi-s}}\leq r_2,\]
    we may use Lemma~\ref{lem-Tmass} applied to \(\nu_y\) and
    \(w\) to find \(w'\in E_y\cap T_y^+ (w,r)\), and in particular lying
    in \(V(y,w-y,r/d_-)\), for which
    \[|w-w'|>\dZero r^{2+\xi-s}=\dZero \dTwo^{2+\xi-s}\psi^{-\frac{1}{2}}\rho.\]
    Hence, as \(|w-u|\leq\rho\),
    \[|w'-u|>\dZero \dTwo^{2+\xi-s}\psi^{-\frac{1}{2}}\rho-\rho
        =(\dZero \dTwo^{2+\xi-s}\psi^{-\frac{1}{2}}-1)\rho\geq \alpha_0 \psi^{-1/2}\rho .\]
%    since \(\psi^{\frac{1}{2}}\leq \sqrt{1-1/\sqrt{2}}\leq \dZero\dTwo^{2+\xi-s}-\alpha_0\).
%
    But
    \[\rho\leq \dOne\psi^{\frac{1}{2}\left(\frac{1}{s-1-\xi}\right)}\leq d_-/\alpha_0< d_-/4\]
    and
    \[r=\dTwo(\psi^{-\frac{1}{2}}\rho)^{\frac{1}{2+\xi-s}}\rho^{-1}\rho
    =\dTwo(\psi^{-\frac{1}{2}}\rho^{s-1-\xi})^{\frac{1}{2+\xi-s}}\rho
    \leq \dTwo\dOne^{\frac{1}{2+\xi-s}}\rho\leq \rho.\]
    Hence, by observation~\ref{obs-two},
    \(w'\in V(y,g,4\rho/d_-)\).

    Similarly, there is \(u'\in E_x\cap T_x^+ (u,r)\) for which
    \[|u'-u|\geq \alpha_0\psi^{-1/2}\rho\]
    and, clearly, \(u'\in V(x,u-x,r/d_-)\).

    Now both \(|u-x|\) and \( |u-y|\) are no less than \( d_- \)
    and
    \[\alpha_0\psi^{-1/2}\rho\leq \alpha_0 d_1 \psi^{\frac{1}{2}\left(
    \frac{2+\xi-s}{s-1-\xi}\right)}<\alpha_0 d_1 \leq d_-,\]
    hence
    \[\min\{|u-x|,\, |u-y|\}>\alpha_0\psi^{-1/2}\rho .\]
    Moreover
    \[r\leq\rho\leq \dOne\psi^{\frac{1}{2}\left(\frac{1}{s-1-\xi}\right)}
    \leq d_1 \psi^{\frac{1}{2}}<\tfrac{1}{20}d_-\psi^{1/2},\]
    and so it follows from observation~\ref{obs-three} that
    \[\emptyset\not= [x,u']\cap [y,w']\subseteq B(u,\alpha_0\psi^{-1/2}\rho).\]
    \begin{figure}
        \begin{center}
            \includegraphics[width=0.7\textwidth, clip=true, keepaspectratio=true]{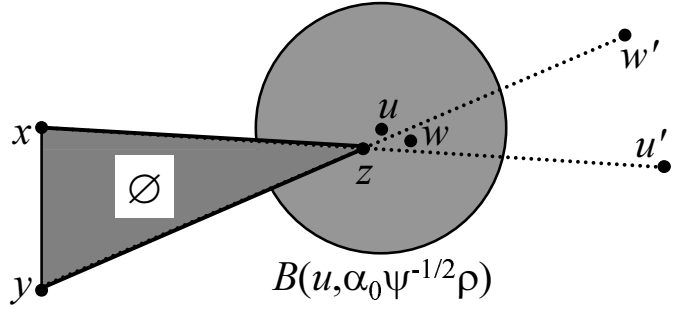}
        \end{center}
        \caption{\(\emptyset\not= [x,u']\cap [y,w']\subseteq B(u,\alpha_0\psi^{-1/2}\rho)\).}
    \end{figure}
    Let \(z\) denote this intersection point.
    Then
    \[([x,z]\cup [x,y]\cup [y,z])\cap\Gamma=\emptyset,\]
    since \(\Gamma\) is connected, \(u'\) is visible from \(x\)
    and \(w'\) is visible from \(y\).
    The observation follows.
\end{oproof}

We now use this observation to find an empty cone with base point
near to \(u\).

\begin{obs}
    \label{obs-five}
    If
    \(0<\rho\leq
    \dOne\psi^{\frac{1}{2}\left(\frac{1}{s-1-\xi}\right)},\) then there is
    \(p\in [\tfrac{1}{2}(x+y),u]\cap B(u,\alpha_0 \psi^{-\frac{1}{2}}\rho)\)
    for which
    \[\mathcal{T}(x,y,p)\cap\Gamma=\emptyset\]
    and
    \[V(p,\tfrac{1}{2}(x+y)-u,\tfrac{2}{5}\psi^{1/2})
    \cap H[x,y;u]\cap\Gamma=\emptyset.\]
\end{obs}
\begin{oproof}[\ref{obs-five}]
    Let
    \[p\in [\tfrac{1}{2}(x+y),u]\cap (V(x,e,4\rho/d_-)\cup
    V(y,g,4\rho/d_-))\]
    be chosen to be at the minimum possible distance from
    \(\tfrac{1}{2}(x+y)\), see Figure~\ref{fig-obs5}.
    Then there is \(\lambda>0\) such that
    \[p=u-\lambda
    \frac{|u-x|e+|u-y|g}{||u-x|e+|u-y|g|}\]
    \begin{figure}
        \begin{center}
            \includegraphics[width=0.7\textwidth, clip=true, keepaspectratio=true]{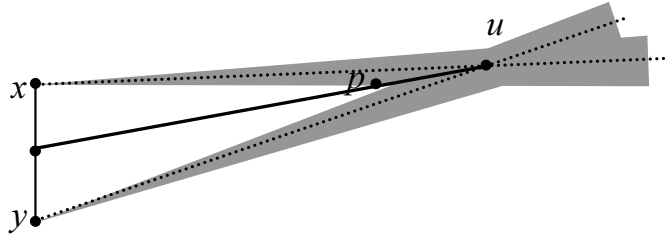}
        \end{center}
        \caption{\(p\in [\tfrac{1}{2}(x+y),u]\cap (V(x,e,4\rho/d_-)\cup
    V(y,g,4\rho/d_-))\).}\label{fig-obs5}
    \end{figure}
    Suppose (without loss of generality) that \(p\in V (y,g,4\rho
    /d_-)\), then there is \(\mu>0\) and \(\sigma\in\{+1,-1\}\)
    such that
    \[p=y+\mu (g+4\rho d_{-}^{-1}\sigma g^\perp).\]
    Hence, if we set \(h=|u-x|e+|u-y|g\), then
    \[u-\lambda h/|h|=y+\mu (g+4\rho d_{-}^{-1}\sigma g^\perp),\]
    which, as \(u-y=|u-y|g\), rearranges to give
    \begin{equation}
        \label{eqn-cl5a}
        |u-y|g-\lambda h/|h|=\mu (g+4\rho d_{-}^{-1}\sigma
        g^\perp).
    \end{equation}
    So taking the inner product of~(\ref{eqn-cl5a}) with
    \(g^\perp\) gives
    \begin{equation}
        \label{eqn-cl5b}
        -\lambda
        \frac{|u-x|}{|h|}\ip{e}{g^\perp}=4\frac{\rho}{d_-}\sigma
        \mu
    \end{equation}
    and taking the inner product of~(\ref{eqn-cl5a}) with \(g\)
    and rearranging gives
    \[|u-y|-\lambda\frac{|u-x|\ip{e}{g}+|u-y|}{|h|}=\mu .\]
    Substituting for \(\mu\) from~(\ref{eqn-cl5b}) gives
    \[|u-y|-\lambda\frac{|u-x|\ip{e}{g}+|u-y|}{|h|}=-\frac{d_-}{4\rho\sigma}
    \cdot\frac{|u-x|\ip{e}{g^\perp}}{|h|}\lambda\]
    and this rearranges to give
    \begin{align*}
    \lambda\left[\left(\frac{d_-}{4\rho\sigma}
    \ip{e}{g^\perp}-\ip{e}{g}\right)|u-x|-|u-y|\right]&=-|u-y|
    |h|\\
    \intertext{and so, substituting for \(h\),}\\
    \lambda\left[(d_-\ip{e}{g^\perp}-4\rho\sigma\ip{e}{g})
    -4\rho\sigma \frac{|u-y|}{|u-x|}\right] &=-4\rho\sigma |u-y| \left|e+\frac{|u-y|}{|u-x|}g\right|
    \end{align*}
    As \( |x-y|\leq d_-/2\), it easily follows that
    \[\frac{2}{3}\leq \frac{|u-x|}{|u-y|},\,\frac{|u-y|}{|u-x|}\leq\frac{3}{2},\]
    and so
    \[
        |\lambda | \leq \frac{4\rho}{d_-}\times\left(1+\tfrac{3}{2}\right)\frac{|u-y|}{|\sigma \ip{e}{g^\perp}
        -4\rho d_{-}^{-1}(\ip{e}{g}+|u-y|/|u-x|)|.}
    \]
    Now \(|\ip{e}{g}+|u-y|/|u-x||\leq 5/2\) and by~(\ref{eqn-cl14}),
    \[|\ip{e}{g^\perp}|=|\ip{e^\perp}{g}|\geq \sqrt{\psi}.\]
    Thus
    \[|\sigma \ip{e}{g^\perp}
        -4\rho d_{-}^{-1}(\ip{e}{g}+|u-y|/|u-x|)| \geq \sqrt{\psi}-10\rho/d_-
        \geq \tfrac{1}{2}\sqrt{\psi},\]
    since \(\rho\leq \dOne\psi^{\frac{1}{2}\left(\frac{1}{s-1-\xi}\right)}
    \leq   d_1\psi^{\frac{1}{2}}\leq (d_-/\alpha_0)\psi^{\frac{1}{2}} < \tfrac{1}{20}d_-\sqrt{\psi}\).
    Hence
    \[|\lambda|\leq 10\rho d_{-}^{-1}|u-y|\times 2\psi^{-1/2}\leq 20(d_+/d_-)
    \psi^{-1/2}\rho \leq \alpha_0\psi^{-\frac{1}{2}}\rho,\]
    and \(p\in B(u,\alpha_0\psi^{-1/2}\rho)\), as claimed.

    Since the hypotheses of observation~\ref{obs-four} are satisfied, there is a point \(z\)
    satisfying its conclusions, and we note that
    \(p\in\mathcal{T}(x,y,z)\). Hence \(\mathcal{T}(x,y,p)\cap\Gamma=\emptyset \).

    To show that
    \[V(p,\tfrac{1}{2}(x+y)-u,\tfrac{2}{5}\psi^{1/2})
    \cap H[x,y;u]\cap\Gamma=\emptyset\]
    we recall that \(h=|u-x|e+|u-y|g\)
    and so
    \[h^\perp =e^\perp |u-x|+g^\perp |u-y|.\]
    We shall estimate \(|\ip{x-u}{h^\perp}|, |\ip{y-u}{h^\perp}|,
    |\ip{x-u}{h}|\) and \(|\ip{y-u}{h}|\) for if \(q\in H(x,y;u)\) satisfies
    \[|\ip{q-p}{h^\perp}|\leq m\ip{q-p}{h}\]
    where \(m=\min\{|\ip{x-u}{h^\perp}/\ip{x-u}{h}|,\, |\ip{y-u}{h^\perp}/
    \ip{y-u}{h}|\}\), then \(q\in\mathcal{T}(x,y,p)\).
    Now
    \begin{align*}
        \ip{x-u}{h^\perp} &= -|x-u| |y-u|\ip{e}{g^\perp}\\
        \ip{x-u}{h} &= -|x-u|^2-|x-u||y-u|\ip{e}{g}\\
        \ip{y-u}{h^\perp} &= -|x-u| |y-u|\ip{g}{e^\perp}\\
        \ip{y-u}{h} &= -|y-u|^2-|x-u||y-u|\ip{e}{g},
    \end{align*}
    and so
    \[\left|\frac{\ip{x-u}{h^\perp}}{\ip{x-u}{h}}\right|=
    \left|\frac{\ip{e}{g^\perp}}{\ip{e}{g}+|x-u|/|y-u|}\right|\]
    and
    \[\left|\frac{\ip{y-u}{h^\perp}}{\ip{y-u}{h}}\right|=
    \left|\frac{\ip{e}{g^\perp}}{\ip{e}{g}+|y-u|/|x-u|}\right|.\]
    Hence
    \(m\geq \tfrac{2}{5}\psi^{1/2}\)
    and the observation follows.
\end{oproof}

We now reach the main part of the proof of the proposition. The
existence of a large empty cone near to \(u\) and \(w\) forces all
other points of \(F_y\cap B(u,\rho)\) to lie in a narrow strip in
direction \(w-y\).

\begin{obs}
    \label{obs-six}
    Let
    \(0<\rho\leq \dOne \psi^{\frac{1}{2}\left(\frac{1}{s-1-\xi}\right)}\)
    and \(r=\dTwo(\psi^{-\frac{1}{2}}\rho)^{\frac{1}{2+\xi-s}}\).
    If \(v,w\in F_y\cap B(u,\rho)\), then \(v\in V(y,w-y,3r/d_-)\).
\end{obs}
\begin{oproof}[\ref{obs-six}]
    Suppose that \(\ip{v-y}{(w-y)^\perp}>0\). (If not, then interchange \(v\) and
    \(w\) --- note that \( \ip{v-y}{(w-y)^\perp}\not=0\) since \(v\) and \(w\) are
    both visible from \(y\).) By observation~\ref{obs-four}, there
    is
    \[z\in V(x,e,r/d_-)\cap V(y,f,r/d_-)\cap B(u,\alpha_0 \psi^{-1/2}\rho)\]
    for which \(\Gamma\cap\mathcal{T}(x,y,z)=\emptyset\).

    By lemma~\ref{lem-Tmass}, we can find \(v'\in E_y\cap T_y^- (v,r)\)
    for which
    \begin{equation}
        \label{eqn-cl6c}
        |v'-u|\geq |v'-v|-|u-v| >\dZero r^{2+\xi-s}-\rho =
        (\dZero\dTwo^{2+\xi-s}\psi^{-\frac{1}{2}}-1)\rho\geq
        \alpha_0\psi^{-1/2}\rho.
    \end{equation}
    We  show that if
    \begin{equation}
        \ip{v-y}{(w-y)^\perp}\geq 3(r/d_-) \ip{v-y}{w-y},\label{eqn-cl6d}
    \end{equation}
    then \(v'\in \mathcal{T}(x,y,z)\), which is impossible.
    To do this, it is enough to show
    \begin{equation}
        \label{eqn-cl6a}
        \ip{v'-x}{e^\perp} <-(r/d_-) \ip{v'-x}{e}
    \end{equation}
    and
    \begin{equation}
        \label{eqn-cl6b}
        \ip{v'-y}{(w-y)^\perp}>(r/d_-) \ip{v'-y}{w-y},
    \end{equation}
    since \(z\in V(x,e,r/d_-)\cap V(y,f,r/d_-)\).

    For~(\ref{eqn-cl6a}): If
    \[v'\in V(y,v-y,r/d_-)\cap V(x,e,r/d_-),\]
    then observation~\ref{obs-three} applied to \(v\) implies that
    \(|v'-u|<\alpha_0\psi^{-1/2}\rho\)
    contradicting~(\ref{eqn-cl6c}). Hence, as \(v'\in
    V(y,v-y,r/d_-)\), we deduce
    \[|\ip{v'-x}{e^\perp}|>(r/d_-) \ip{v'-x}{e}\]
    and it only remains to show that \(\ip{v'-x}{e^\perp}<0\).
    If \(\ip{v'-x}{e^\perp}\geq 0\), then
    \( \ip{v'-x}{e^\perp}> (r/d_-)\ip{v'-x}{e}\). But, by
    observation~\ref{obs-two} applied to \(v\),
    \(|\ip{v'-y}{g^\perp}|\leq 4(\rho/d_-) \ip{v'-y}{g}\) and so, as
    \(|v'-y|<|v-y|\), we find
    \begin{align*}
        |\ip{v'-y}{g^\perp}| <4\rho/d_- |v-y| &\leq 4\rho/d_-
        (|u-y|+\rho )\\
        &\leq (5\rho/d_-) |u-y|,\mbox{ as \(\rho\leq d_-/4\)}\\
        &\leq 5 (d_+/d_-)\rho
    \end{align*}
    and
    \( \ip{v'-u}{g^\perp}=\ip{v'-y}{g^\perp}\). So
    \[|\ip{v'-u}{g^\perp}|\leq 5 ( d_+/d_-)\rho .\]
    Now
    \[\ip{v'-u}{g}=\ip{v'-y}{g}-|u-y|\leq |v-y|-|u-y|\leq\rho .\]
    Let \(q\) be the point of intersection of \([x,u]\) with
    \([y,y+|y-u|(g+4(\rho/d_-)g^\perp)]\).
    Then since \(v'\in V(y,g,4(\rho/d_-))\) and \(\ip{v'-x}{e^\perp}>
    (r/d_-)\ip{v'-x}{e}\), it follows that
    \(\ip{v'-u}{g}\geq \ip{q-u}{g}\). Hence it is enough to find a lower bound for \(\ip{q-u}{g}\).
    There is \(0<\lambda
    <|u-y|\) and \(0<\mu<|x-u|\) such that
    \[q=y+\lambda (g+4(\rho/d_-) g^\perp)=x+\mu e,\]
    and so
    \[-|u-y|g+\lambda(g+4(\rho/d_-) g^\perp)=(\mu-|u-x|)e.\]
    Taking inner products of this expression with \(g\) and \(g^\perp\), and
    solving for \(\lambda\) gives
    \[\lambda =|u-y| \left(1+4\left(\frac{\rho}{d_-}\right)
    \frac{\ip{e}{g}}{\ip{e^\perp}{g}}\right)^{-1}.\]
    Hence
    \[\lambda \geq |u-y| (1+4(\rho /d_-)\psi^{-1/2})^{-1}
    \geq |u-y| (1-4(\rho /d_-)\psi^{-1/2}).\]
    Thus
    \[\ip{q-u}{g}\geq -4(\rho /d_-)\psi^{-1/2}|u-y|\]
    and so
    \[\ip{v'-u}{g}\geq -4(\rho /d_-)\psi^{-1/2}|u-y|\geq
    -4(d_+/d_-)\psi^{-\frac{1}{2}}\rho.\]
    Hence
    \[|v'-u|\leq 5\rho\tfrac{d_+}{d_-}+4\rho\psi^{-1/2}\tfrac{d_+}{d_-}\leq
    9(d_+/d_-)\psi^{-1/2}\rho<\alpha_0\psi^{-1/2}\rho,\]
    a contradiction.

    For~(\ref{eqn-cl6b}), notice that \(v'-y=\alpha (v-y)+\beta (v-y)^\perp\)
    for some \(0<\alpha <1\) and \(|\beta| < \alpha ( r/d_-)\).
    Thus, using~(\ref{eqn-cl6d}),
    \begin{align*}
        \ip{v'-y}{(w-y)^\perp} &= \alpha \ip{v-y}{(w-y)^\perp}+
        \beta \ip{(v-y)^\perp}{(w-y)^\perp}\\
        &\geq 3\alpha (r/d_-) \ip{v-y}{w-y}- |\beta||\ip{v-y}{w-y}|\\
        &\geq 2\alpha r/d_- \ip{v-y}{w-y}
    \end{align*}
    and
    \[\ip{v'-y}{w-y}=\alpha \ip{v-y}{w-y}+\beta\ip{(v-y)^\perp}{w-y}.\]
    But
    \[\ip{(v-y)^\perp}{w-y} =\ip{(v-y)^\perp}{w-v}\]
    and so \(|\ip{(v-y)^\perp}{w-y}|\leq 2\rho |v-y|\), and
    \[\ip{v-y}{w-y}=\ip{v-y}{w-v}+|v-y|^2\]
    and so
    \[|\ip{v-y}{w-y}|\geq |v-y| (|v-y|-2\rho)\geq (d_-/2)|v-y|,\]
    since \(|v-y|\geq d_-\) and \(\rho\leq d_-/4\).
    Thus
    \[|\ip{(v-y)^\perp}{w-y}|\leq 2\rho |v-y|
    \leq 4(\rho/d_-) \ip{v-y}{w-y} .\]
    So
    \begin{align*}
        \ip{v'-y}{w-y} &\leq (\alpha +4(\rho/d_-)
        |\beta|)\ip{v-y}{w-y}\\
        &\leq (1+4(\rho/d_-)(r/d_-))\alpha \ip{v-y}{w-y}\\
        &< (1+4(\rho/d_-)(r/d_-)) (d_-/(2r)) \ip{v'-y}{(w-y)^\perp}
    \end{align*}
    and rearranging gives
    \[\ip{v'-y}{(w-y)^\perp}> (2r/d_-)(1+4(\rho/d_-)(r/d_-))^{-1}\ip{v'-y}{w-y}\]
    which, since \(4(\rho/d_-)(r/d_-)\leq r/d_-\leq r_2/d_-\leq 1\), implies~(\ref{eqn-cl6b}).
\end{oproof}
    We can now finish the
    proof of the proposition.

    Let \(0<\rho\leq \dOne \psi^{\frac{1}{2}\left(\frac{1}{s-1-\xi}\right)}\)
    and \(r=\dTwo (\psi^{-\frac{1}{2}}\rho)^{\frac{1}{2+\xi-s}}\).
    Suppose  \(w\in B(u,\rho)\cap F_y\), then
    \[F_y\cap B(u,\rho)\subset V(y, (w-y),3r/d_-).\]
    Thus \(F_y\cap B(u,\rho)\) is contained in a rectangle of
    height \(2\rho\) and width \(6rd_+/d_-\) which can be covered
    by \((2+2\rho/r)(2+6d_+/d_-)\) boxes of side \(r\).
    Since
    \[(2+2\rho/r)(2+6d_+/d_-)\leq 32 (\rho/r)d_+/d_-,\]
    and
    \(\sqrt{2}r\leq \sqrt{2}\alpha_1 \left(d_1\psi^{\frac{1}{2}\left(\frac{2+\xi-s}{s-1-\xi}\right)}\right)^{\frac{1}{2+\xi- s}} \leq \sqrt{2}\dTwo\dOne^{\frac{1}{2+\xi-s}}\leq \sqrt{2}r_2\leq r_0\), we estimate that
    \[\nu_y (F_y\cap B(u,\rho))\leq 2^{5+s/2}(d_+/d_-)\rho
    r^{s-1}=2^{5+s/2}\dTwo^{s-1}(d_+/d_-)\psi^{-\frac{1}{2}\left(\frac{s-1}{2+\xi-s}
    \right)}\rho^{\frac{1+\xi}{2+\xi-s}},\]
    as required. The remainder of the proposition follows from
    observation~\ref{obs-five}.
\end{proof}

\subsection{Mass estimate proposition}

The main utility of Proposition~\ref{prop-geom1} lies in its use
in proving the following proposition.

\begin{prop}
    \label{prop-mass}
    Let \(\Gamma\) be a non-empty compact connected subset of
    \(\R^2\) and let \(A\) and \(B\) be compact subsets of \(\R^2\).
    Suppose that \(s>1\), \(0<\xi<s-1\), \(0<r_1\leq r_0\leq 1\), \(0<d_-\leq d_+\) with \(d_-\leq 1\) and
    \(M>0\) are given.
    Let \(x,y\in \R^2\setminus\Gamma\) satisfy
    \[0<2|x-y|<d_-\leq \min\{d(x,\Gamma),\, d(y,\Gamma)\}
        \leq \max\{ d(x,\Gamma),\, d(y,\Gamma)\}+|\Gamma|\leq d_+ .\]
    Let \(\nu_x\) and \(\nu_y\) be Radon measures supported in \(\Gamma_x\)
    and \(\Gamma_y\) respectively and let
    \[F_x\subseteq E_x\subseteq \Gamma_x\mbox{ and } F_y\subseteq E_y\subseteq
    \Gamma_y\]
    be compact sets.
    Suppose that:
    \begin{enumerate}
        \item for all \(u\in E_x\), \(v\in E_y\) and \(0<r\leq
        r_0\) both
        \[\nu_x (B(u,r)\leq r^s\mbox{ and }\nu_y (B(v,r))\leq r^s;\]
        \item for all \(u\in F_x\), \(v\in F_y\) and \(0<r\leq
        r_1\) both
        \[\nu_x (T_x^\pm(u,r)\cap E_x)\geq Mr^{1+\xi}\mbox{ and }
        \nu_y (T_y^\pm(v,r)\cap E_y)\geq Mr^{1+\xi};\]
        \item there is \(\psi\in (0,1/{2})\) such that for
        \(u\in (F_x\cap A)\cup (F_y\cap B)\),
        \[\ip{(u-x)^\wedge}{(u-y)^\wedge}\in [1/{2},1-\psi].\]
    \end{enumerate}
    Then there are constants \(\dEight,\dFour >0\) such that
    if \(0<\rho\leq \dFour\psi^{\frac{1}{2}\frac{1}{s-1-\xi}}\), then
    \begin{align*}
        \lefteqn{(\nu_x\otimes\nu_y) \left( (F_x\times F_y)\cap(A\times
        B)\cap\{
        (u,v): |u-v|\leq \rho\}\right)}\qquad \qquad\qquad \qquad \\
        &\leq \dEight \arcdiam_{\frac{1}{2}(x+y)} (A\cap F_x\cap B(F_y\cap B,\rho))
        (\psi^{-\frac{1}{2}}\rho)^{\frac{s+\xi}{2+\xi-s}}.
    \end{align*}
\end{prop}

\begin{proof}
    Let \(\dFour=\frac{5}{2^{3/2}\alpha_0}\dOne\) and observe that, as
    \(\dFour\leq \dOne\),
    Proposition~\ref{prop-geom1} implies
    \begin{align*}
        \lefteqn{(\nu_x\otimes\nu_y) \left( (F_x\times F_y)\cap(A\times B)\cap\{
        (u,v): |u-v|\leq \rho\}\right)}\\
        & = \int_{B(F_y\cap B, \rho)}
        \nu_y|_{F_y\cap B} \left(\{v: |u-v|\leq \rho\}
            \right)d\nu_x|_{F_x\cap A} (u)\\
        & = \int_{B(F_y\cap B, \rho)}
        \nu_y\left(F_y\cap B\cap B(u,\rho)\right)
            d\nu_x|_{F_x\cap A} (u)\\
            & \leq \dThree \psi^{-\frac{s-1}{2(2+\xi-s)}}\rho^{\frac{1+\xi}{2+\xi-s}}
            \nu_x \left(F_x\cap A\cap B(F_y\cap B, \rho)
            \right).
    \end{align*}
    It remains to estimate
    \[\nu_x \left(F_x\cap A\cap B(F_y\cap B, \rho)
            \right) .\]
    We begin by noticing that for each
    \[u\in F_x\cap A\cap B(F_y\cap B, \rho),\]
    Proposition~\ref{prop-geom1} guarantees the existence of
    \(p_u\in  [\tfrac{1}{2}(x+y), u]\cap B (u,\alpha_0 \psi^{-1/2}\rho)\)
    such that \(\mathcal{T}(x,y,p_u)\cap \Gamma=\emptyset\)
    and
    \[V(p_u,\tfrac{1}{2}(x+y)-u,\tfrac{2}{5}{\psi^{\frac{1}{2}}})\cap
    \Gamma\cap H(x,y;u)=\emptyset .\]

    Let \(\sigma =\dFive\rho\) and fix \(v\in F_x\cap A\cap B(F_y\cap B,
    \rho)\).
    Then Lemma~\ref{lem-intercone} guarantees that if
    \[w\in V(\tfrac{1}{2}(x+y), p_v-\tfrac{1}{2}(x+y),\sigma)\setminus V(p_v,-(p_v-\tfrac{1}{2}(x+y)),
    \tfrac{2}{5}\psi^{\frac{1}{2}}),\]
    then
    \begin{align*}
        \ip{w-p_v}{(p_v-\tfrac{1}{2}(x+y))^\wedge}&\geq
        -\frac{\sigma}{\sigma+\frac{2}{5}\psi^{\frac{1}{2}}}
        |p_v-\tfrac{1}{2}(x+y)|\\
        &  \geq -\tfrac{5}{2}d_+\left(\dFive\right)\psi^{-\frac{1}{2}}\rho =-\alpha_0\psi^{-\frac{1}{2}}\rho.
    \end{align*}

    So suppose \(u,v\in F_x\cap A\) with \(u\in V(\tfrac{1}{2}(x+y),v,\sigma )\)
    and assume, without loss of generality, that
    \[|u-\tfrac{1}{2}(x+y)|\leq |v-\tfrac{1}{2}(x+y)|.\]
    We wish
    to estimate \(\ip{u-v}{(v-\tfrac{1}{2}(x+y))^\wedge}\)
    from below. (An easy upper bound is given by zero.) From the preceding we know that
    \[\ip{u-p_v}{(v-\tfrac{1}{2}(x+y))^\wedge}=\ip{u-p_v}{(p_v-\tfrac{1}{2}(x+y))^\wedge}
    \geq -\alpha_0\psi^{-\frac{1}{2}}\rho.\]
    Hence
    \begin{align*}
        \ip{u-v}{(v-\tfrac{1}{2}(x+y))^\wedge} &=
            \ip{u-p_v}{(v-\tfrac{1}{2}(x+y))^\wedge}
            +\ip{p_v-v}{(v-\tfrac{1}{2}(x+y))^\wedge}\\
            &\geq
            -\alpha_0\psi^{-\frac{1}{2}}\rho-\alpha_0
            \psi^{-\frac{1}{2}}\rho\\
            &\geq -2\alpha_0\psi^{-\frac{1}{2}}\rho.
    \end{align*}
    Thus
    \[V(\tfrac{1}{2}(x+y),v-\tfrac{1}{2}(x+y),\sigma)\cap (F_x\cap A\cap
        B(F_y\cap B,\rho )\cap B(\tfrac{1}{2}(x+y),|v-\tfrac{1}{2}(x+y)|))\]
    can be covered by
    \[\frac{2\alpha_0\psi^{-1/2}\rho}{\frac{2}{5}\alpha_0 \rho}=5\psi^{-\frac{1}{2}} \]
    boxes of side \(\frac{4}{5}\alpha_0 \rho\).
    Hence, by using the mass estimate in Proposition~\ref{prop-geom1}
    (for \(\nu_x\) rather than \(\nu_y\)), since
    \(2\sqrt{2} \frac{2}{5}\alpha_0 \rho
    \leq \dOne \psi^{\frac{1}{2}\frac{1}{s-1-\xi}}\),
    \[\nu_x \left( V(\tfrac{1}{2}(x+y),v-\tfrac{1}{2}(x+y),\sigma)
    \cap F_x\cap A\cap B(F_y\cap B, \rho)\cap B(\tfrac{1}{2}(x+y),|v-\tfrac{1}{2}(x+y)|) \right)\]
    is at most
    \[5\psi^{-\frac{1}{2}}\times \dThree \psi^{-\frac{1}{2}\frac{s-1}{2+\xi-s}}
    \left(\tfrac{4\sqrt{2}}{5}\alpha_0 \rho\right)^{\frac{1+\xi}{2+\xi-s}}
            =5\dThree\left(\tfrac{4\sqrt{2}}{5}\alpha_0\right)^{\frac{1+\xi}{2+\xi-s}}
            (\psi^{-\frac{1}{2}}\rho)^{\frac{1+\xi}{2+\xi-s}}.\]
    By choosing \(v\) to be as far from \(\tfrac{1}{2}(x+y)\) as possible and counting the number of such cones needed to cover \(F_x\cap A\), we obtain
    \begin{align*}
        \lefteqn{\nu_x (F_x\cap A \cap B(F_y\cap B,\rho))}\\
        &\leq {2\arcdiam_{\frac{1}{2}(x+y)} (F_x\cap A \cap B(F_y\cap
        B,\rho))}\sigma^{-1}
        \times 5\dThree\left(\tfrac{4\sqrt{2}}{5}\alpha_0\right)^{\frac{1+\xi}{2+\xi-s}}
            (\psi^{-\frac{1}{2}}\rho)^{\frac{1+\xi}{2+\xi-s}}\\
        &=\dEight
    \arcdiam_{\frac{1}{2}(x+y)} (F_x\cap A \cap B(F_y\cap B,\rho))(\psi^{-\frac{1}{2}}\rho)^{\frac{1+\xi}{2+\xi-s}}\rho^{-1}
    \end{align*}
    for  \(\dEight=25\dThree d_+(4\sqrt{2}/5)^{\frac{1+\xi}{2+\xi-s}}\alpha_0^{\frac{s-1}{2+\xi-s}}\), and this implies the claim.
\end{proof}

\section{Measurability results}\label{sec-meas}

In this section we prove the measurability of various maps that we
use in the proof of Theorem~\ref{thmresult}. In particular, we
show that if \(B\) is a compact set that is disjoint from
\(\Gamma\), then there is a universally-measurable map that
assigns to each point \(x\in B\) for which \(\Gamma_x\) has large
dimension, a Radon measure of large dimension that is `supported'
on \(\Gamma_x\).

Let \(B\) be a compact subset of the plane disjoint from the
non-empty compact connected set \(\Gamma\). Letting \(S^1\) denote
the unit circle, we define \(K\subseteq B\times S^1\times \R^+\)
by
\[K=\{ (x,\theta,t )\in B\times S^1\times\R^+: x+t\hat{\theta}\in\Gamma \}\]
where \(\hat{\theta}=(\cos\theta, \sin\theta)\). Notice that \(K\)
is  compact and a lifting of \(\Gamma\).

For \(x\in B\) and \(\theta\in S^1\), define \(\gamma\colon
B\times S^1\to \R^+\cup\{\infty\}\) by
\[\gamma (x,\theta)=\begin{cases}
\inf\{ t>0: x+t\hat{\theta} \in\Gamma\} & \mbox{if }(x+\R
\hat{\theta})\cap\Gamma\not=\emptyset\\
\infty &\mbox{otherwise.}
\end{cases}.\]
  Observe that \( x+\gamma
(x,\theta)\hat{\theta}\in \Gamma\) for any \(x\in B\) and
\(\theta\in S^1\) for which \( (x+\R\hat{\theta})
\cap\Gamma\not=\emptyset\).

Let
\[\gr(\gamma)=\{(x,\theta,\gamma(x,\theta)): (x,\theta )\in B\times
S^1, \gamma (x,\theta)<\infty\},\] then
\[\gr(\gamma)\subseteq
    K \subseteq
    B\times S^1\times I,\] where \(I=[0,\diam
(B\cup\Gamma)] \).
\begin{lem}\label{lem-grgamma}
    The function \(\gamma\) is lower semi-continuous. In
    particular, \(\gr(\gamma)\) is a \(G_\delta\)-subset of \(K\).
\end{lem}

\begin{proof}
    That \(\gamma\) is lower semi-continuous follows
    readily from the observation that its graph is the lower
    envelope of the compact set \(K\).

    The fact that \(\gr(\gamma)\) is \(G_\delta\) is a standard
    result concerning functions of Baire class 1, see for
    example,~\cite[Ch~II,\S 31 VII, Theorem~1]{kuratowski}.
\end{proof}

For \(C\subseteq K\) and \(x\in B\), let \(C_x\) be given by
\[C_x = C\cap \left(\{x\}\times S^1\times I\right),\]
 the slice of \(C\) through \(x\). For ease, we let \(\gr_x (\gamma)\)
 denote \( (\gr (\gamma))_x\).

 Recall that \(\mathcal{M}(K)\) denotes the Radon measures
supported in \(K\). The set \(\mathcal{M}(K)\) can be given the
topology of weak convergence by using as a base, sets of the form
\[\left\{\mu\in \mathcal{M} (K): \int f \, d\mu <a\right\},\]
where \(a\in\R\) and \(f\in C(K)\), the set of real-valued
continuous functions on \(K\).
 It
turns out that \(\mathcal{M}(K)\) with this topology is a Polish
space, see~\cite[\S 14.15]{mat1} and \cite[II.17]{kechris}.

\begin{lem}\label{prop-meas1}
    Let \(E\) be a Borel subset of \(K\).
    Then the functions \(F_E\colon \mathcal{M}(K)\to\R\) and
    \(G_E\colon B\times\mathcal{M}({K})\to \R\)
    given by
    \[F(\nu)=\nu (E)\mbox{ and }G_E (x,\nu)=\nu (E_x)\]
    are Borel.

    In particular,
    \[\{ (x,\nu)\in B\times \mathcal{M}({K}) : \nu (E_x) >0 \}\]
    is a Borel set.
\end{lem}

\begin{proof}
    Let \(E\subseteq K\) be a Borel set.
    We show that \(G_E\) is Borel; the proof
    that \(F_E\) is Borel is similar.

    Suppose first that \(E\) is a compact subset of \(K\).
    Then for \(x\in B\), \(E_x\) is also compact, and
    for \(\mu\in \mathcal{M}(K)\),
    \begin{align*}
        \lefteqn{G_E(x,\nu )=\nu (E_x)<c \mbox{ if and only if}}&\\
        &\qquad\mbox{ there is }f\in C^+(K)\mbox{ such that }
        f>1\mbox{ on }E_x\mbox{ and
        }\int\! f\, d\nu <c.
    \end{align*}
(Here \(C^+ (K)\) denotes the set of non-negative real-valued
continuous functions on \(K\).)
    For a given \(f\in C^+(K)\), the sets
    \[B_f=\{x\in B : f>1\mbox{ on }C_x\}\]
    and
    \[M_f=\{ \nu \in \mathcal{M}(K): \int\! f\, d\nu <c\}\]
    are open subsets of \(B\) and \(\mathcal{M}(K)\),
    respectively. Hence
    \[\{(x,\nu):\nu (E_x)<c\}=\bigcup_{f\in C^+(K)}B_f\times M_f\]
    is an open set, and so \(G_E\) is upper semi-continuous, and in
    particular, Borel.

    If \(E_1\subseteq E_2\subseteq E_3\subseteq\cdots\) is an
    increasing
    sequence of compact sets, and
    \(G_1,G_2,G_3,\ldots\), the associated sequence of maps, then
    \[G_{\cup_{i\in\N} E_i}=\lim_{i\to\infty} G_i\]
    is a Borel map. Similarly, if \(E_1\supseteq E_2\supseteq E_3\supseteq\cdots\) is a
    decreasing
    sequence of compact sets, and
    \(G_1,G_2,G_3,\ldots\), the associated sequence of maps, then
    \[G_{\cap_{i\in\N} E_i}=\lim_{i\to\infty} G_i\]
    is also a Borel map. It follows that for a general Borel set
    \(E\), the map \(G_E\) is Borel, as required.
\end{proof}

\begin{lem}\label{lem-meas4}
    Let \(E\) be a Borel subset of \(K\).
    The set
    \[\{(x,\nu)\in B\times\mathcal{M}(K): \nu (K\setminus (\gr_x (\gamma)\cap E))=0 \}\]
    is Borel in \(B\times\mathcal{M}(K)\).
\end{lem}

\begin{proof}
    Observe that
    \begin{align*}
        \lefteqn{\{(x,\nu)\in B\times\mathcal{M}(K): \nu (K\setminus (\gr_x
        (\gamma)\cap E))=0 \}} \\
        &= \{ (x,\nu): \nu ((\gr_x
        (\gamma)\cap E))=\nu (K)\}\\
        &= \{(x,\nu): \nu ((\gr
        (\gamma)\cap E)_x)-\nu (K)=0\}\\
        & = \{ (x,\nu): G_{\gr
        (\gamma)\cap E} (x,\nu)-F_K (\nu)=0\}.
    \end{align*}
    Hence, since lemma~\ref{prop-meas1} implies \(F_K\) and \(G_{\gr
        (\gamma)\cap E}\) are Borel functions,
    this set  is Borel.
\end{proof}

Define \(\Pi\colon B\times S^1\times I\to \R^2\) by
\[\Pi (x,\theta, t)=x+t\hat{\theta}\]
and observe that \(\Pi\) is continuous.

\begin{lem}
    For \(x\in B\), if \(A\subseteq (\{x\}\times S^1\times I)\cap K\), then
    \[\dimh (\Pi (A))=\dimh (A).\]
\end{lem}

\begin{proof}
    This follows from the fact that \(\Pi\) is bi-Lipschitz when restricted to \(\{x\}\times S^1\times I\).
\end{proof}

In particular, since
\[\Gamma_x=\Pi ( (\{x\}\times S^2\times I)\cap \gr (\gamma))
=\Pi (\gr_x (\gamma)),\] it follows that
\[\dimh (\Gamma_x)
    =\dimh (\gr_x (\gamma)),\]
for each \(x\in B\). Recall that for \(A\subseteq K\) and
\(s\in\R\),
\begin{align*}
    \lefteqn{\mathcal{M}^s (A)}\\
    &=\{\nu\in\mathcal{M}(K) : \nu (A)>0 \mbox{ and }
    \nu (B(\zeta ,r))\leq r^s \mbox{, for }\zeta\in K ,\, r\in (0,1] \}.
\end{align*}
It is an easy calculation, which we omit, to check that
\(\mathcal{M}^s(K)\) is a Borel set. Since \(\gr_x (\gamma)\) is a
Borel set,
    \begin{align*}
        \dimh (\gr_x (\gamma))&=\sup\{\sigma: \mathcal{M}^{\sigma}(\gr_x (\gamma))
        \not=\emptyset\}.
    \end{align*}

\begin{prop}\label{cor-analytic}
    Let \(C\) be a Borel subset of the plane. Then for \(t\geq 0\),
    \[\{x\in B :\dimh (\Gamma_x\cap C)>t\}\]
    is an analytic set.
\end{prop}

\begin{proof}
    Let \(E\subseteq K\) be given by
    \(E=\Pi^{-1}(C)\cap
    K\),
    and observe that \(\gr (\gamma)\cap E\) is a Borel subset of \(K\).
    For \(t\geq 0\)
    \begin{align*}
        \{x:\dimh (\Gamma_x\cap C)>t\} &=\{x :\dimh ((\gr (\gamma)\cap E)_x)>t\}
        \\&= \{x: \mathcal{M}^\tau ((\gr (\gamma)\cap E)_x)\not=\emptyset
        \mbox{ for some }\tau >t \}\\
        &= \bigcup_{p\in\Q^+}\{x: \mathcal{M}^{t+p} ((\gr (\gamma)\cap E)_x)\not=\emptyset \}.
    \end{align*}
    However, if \(\pi_B\colon B\times \mathcal{M}(K)\to B\) denotes
    the coordinate projection onto \(B\), then
    \begin{align*}
    \lefteqn{\{x: \mathcal{M}^{t+p} ((\gr (\gamma)\cap E)_x)\not=\emptyset \} }&\\
    &=\pi_B (\{ (x,\nu): x\in B,\, \nu\in
    \mathcal{M}^{t+p}((\gr (\gamma)\cap E)_x)\} )\\
    &=\pi_B\left(\{ (x,\nu)\in B\times\mathcal{M}^{t+p}(K):
        \nu ((\gr (\gamma)\cap E)_x)>0\})\right).
    \end{align*}
    Hence lemmas~\ref{lem-grgamma} and~\ref{prop-meas1} together imply
    that
    \(\{x: \mathcal{M}^{t+p} ((\gr (\gamma)\cap E)_x)\not=\emptyset \}\)
    is the coordinate-wise projection of a Borel set from a product of  Polish
    spaces, and so it is
    analytic, see~\cite[Chapter~III]{kechris}. Hence
    \(\{x:\dimh (\Gamma_x\cap C)>t\}\) is also analytic.
\end{proof}

Our last result in this section is a selection theorem and allows
us to choose, in a measurable way, an element of \(\mathcal{M}^t
(\gr_x (\gamma))\) whenever \(x\in B\) is such that \(\dimh
(\Gamma_x)>t\).

\begin{prop}\label{prop-unif}
    Let \(t\geq 0\) and \(C\) be a Borel subset of the
    plane.
    There is a map
    \begin{align*}
        \omega\colon \{x\in B:\dimh (\Gamma_x\cap C)>t\} &\to
        \mathcal{M}(K)\\
        x &\mapsto \omega_x
    \end{align*}
    such that:
    \begin{enumerate}
        \item \(\omega\) is \(\sigma(\bof{A})\)-measurable,
        \item \(\omega_x\in \mathcal{M}^t (\gr_x (\gamma)\cap \Pi^{-1}(C))\) for each \(x\), and
        \item \(\omega_x (K\setminus (\gr_x (\gamma)\cup \Pi^{-1}(C)))=0\) for each \(x\).
    \end{enumerate}
    (Here \(\sigma(\bof{A})\) denotes the
    \(\sigma\)-algebra generated by the analytic sets in \(B\).)

    In particular, \(\omega\) is $\mu$-measurable for every Radon
    measure \(\mu\) on \(B\).
\end{prop}
\begin{proof}
     Let \(E=\Pi^{-1}(C)\cap K\), a Borel set.
    Since
    \[(B\times \mathcal{\sigma}(K))\cap\{(x,\nu):\nu ((\gr(\gamma)\cap E)_x)>0\}
    \cap\{(x,\nu): \nu (K\setminus (\gr(\gamma)\cap E)_x)=0\}\]
    is Borel in \(B\times\mathcal{M}(K)\), claims
    1,2, and 3 follow readily from the Jankov-von Neumann Uniformisation
    Theorem. (See~\cite[Theorem~18.1]{kechris} for a statement of
    this theorem.)

    See~\cite[Theorem~21.10]{kechris} for a proof of Lusin's
    Theorem that analytic sets are universally measurable, from
    which it follows that sets in the \(\sigma\)-algebra generated
    by analytic sets are also universally measurable.
\end{proof}

\section{Proof of Theorem~\ref{thmresult}}
We now draw our preparatory work together and prove
Theorem~\ref{thmresult}.

Let \(\Gamma\) be a compact connected subset of the plane for
which \(1<\dimh(\Gamma)\leq 2\). If \(\dimh (\Gamma )=2\), then
let \(d=2\), otherwise choose \(\dimh (\Gamma)<d<2\). Notice that
in both cases this implies that whenever \(\nu\) is a non-zero
Radon measure supported in \(\Gamma\), then
\[I_d (\nu)=+\infty.\]
(If \(d=2\), then, since \(\mathcal{H}^2 (\Gamma)<\infty\),
Theorem~8.7 of~\cite{mat1} implies \( I_2 (\nu)=+\infty\).)

\subsection{Measure theoretic decomposition}\label{ssec-decomp}

Fix \(d>s>1\) and let
\(\emptyset\not=B^{(0)}\subseteq\R^2\setminus\Gamma\) be a compact
set for which \(\diam (B^{(0)}) \leq \frac{1}{100}\dist
(B^{(0)},\Gamma )\). It is enough for us to show that
\[\dimh(\{x\in B^{(0)}: \dimh (\Gamma_x)>s\})<\tfrac{1}{2}+\sqrt{d-\tfrac{3}{4}}.\]

%From Proposition~\ref{cor-analytic}, we see that
%\[E=\{x\in B: \dimh (\Gamma_x)>s\}\]
%is an analytic set.

Since \(\Gamma\) is compact, we can find finitely many open sets
\(U_1\), \(U_2\),\ldots, \(U_N\) that intersect \(\Gamma\) such
that \(\Gamma\subseteq \bigcup_{i=1}^N U_{i}\) and \(\diam
(U_i)\leq \tfrac{1}{100} \dist (B^{(0)},\Gamma )\) for each \(i\).

It follows that
\[E=\bigcup_{i=1}^N\{x\in B^{(0)}: \dimh (\Gamma_x\cap U_i)>s\}=\bigcup_{i=1}^N E_i\mbox{, say.}\]
Clearly each \(E_i\) satisfies
\[\diam (E_i)\leq \tfrac{1}{100} \dist (E_i,\Gamma) \leq \tfrac{1}{100} \dist (B^{(0)},\Gamma).\]
From Proposition~\ref{cor-analytic}, we see that each \(E_i\) is
an analytic set.

Moreover, for \(t>0\), if
    \[\dimh \left(\{x\in B^{(0)}: \dimh (\Gamma_x)>s\}\right)>t,\]
then we can find an \(i\)  such that
\begin{equation}
    \label{eqn-decomp1}
    \dimh(E_{i})>t .
\end{equation}
So suppose \(t>0\) and \(i\) are such that \( \dimh(E_{i})>t\).
Our objective is to find an upper bound for the size of \(t\) in
terms of \(d\) and \(s\).

By Theorem~\ref{thmhdim}, there is a nonzero Radon measure \(\mu\)
with compact support \(B^{(1)}\subseteq A_i\subseteq B^{(0)}\)
such that whenever \(x\in \R^2\) and \(r>0\), then
\begin{equation}
    \label{eqn-muBxr}
    \mu (B(x,r))\leq r^t .
\end{equation}

Proposition~\ref{prop-unif} enables us to find a
\(\sigma(\bof{A})\)-measurable function
\begin{align*}
    \omega\colon B^{(1)} &\to\mathcal{M}^s (K)\\
     x&\mapsto \omega_x,
\end{align*}
(where \(K=\{ (x,\theta,t )\in B^{(1)}\times S^1\times\R^+:
x+t\hat{\theta}\in\Gamma \}\)) such that
\begin{itemize}
    \item \(\omega_x (\gr_x (\gamma)\cap\Pi^{-1}(U_i))>0\),
    \item \(\omega_x (K\setminus (\gr_x (\gamma)\cap \Pi^{-1}(U_i) ))=0\),
    for each \(x\in B^{(1)}\).
\end{itemize}
Moreover, there is a constant \(C\) such that \(\omega_x (K)\leq
C\) for all \(x\). By Lusin's Theorem~\cite[2.3.5]{federer}, there
is a compact set \(B^{(2)}\subseteq B^{(1)}\) such that
\begin{itemize}
    \item \(\mu (B^{(2)} )>0\), and
    \item \(\omega |_{B^{(2)}}\) is a continuous map.
\end{itemize}
Let \(\mu^{(2)}=\mu |_{B^{(2)}}\) and for Borel \(E\subseteq K\)
define
\[m^* (E)=\int\!\omega_x (E)\, d\mu^{(2)} (x) ,\]
and extend \(m^*\) to arbitrary \(A\subseteq K\) by setting
\(m^*(A)=\inf\{m^* (E): A\subseteq E\mbox{ and }E\mbox{ is
Borel}\}\). We omit the routine verification that
    \(m^*\) is a Radon measure on \(K\).

For \(x\in B^{(2)}\) define a Radon measure \(\nu_x\) on
\(\Gamma\) by
\[\nu_x (A)=\omega_x (\Pi^{-1}(A))=\omega_x (\Pi^{-1}(A)\cap (\{x\}\times S^1\times I))\mbox{, for }A\subseteq \R^2,\]
and observe that the continuity of the map \(\omega\) implies that
\(x\mapsto \nu_x\) is a Borel measurable function. Also notice:
\begin{itemize}
    \item for \(x\in B^{(2)}\), \(\nu_x (\R^2\setminus \Gamma_x)=0\) and \(0<\nu_x (\R^2)\leq C\),
    \item for \(x\in B^{(2)}\), \(u\in\R^2\) and \(0<r\leq 1\),
    \(\nu_x (B(u,r))\leq r^s\).
\end{itemize}

We now analyse the geometry of the measures \(\nu_x\).

Fix \(0<\xi<s-1\). Then for all \(x\in B^{(2)}\),
Lemma~\ref{lem-cones} implies that for \(\nu_x\)-a.e.\ \(u\in
\Gamma_x\),
\begin{equation}
    \label{eqn-tubebd1}
    \min\left\{\liminf_{r\searrow 0}\frac{\nu_x (T^+_x (u,r))}{r^{1+\xi}},\,
    \liminf_{r\searrow 0}\frac{\nu_x (T^-_x
    (u,r))}{r^{1+\xi}}\right\}=+\infty.
\end{equation}
That is, for all \(x\in B^{(2)}\) and \(\omega_x\)-a.e.\
\(\zeta\in K\),
\[\liminf_{r\searrow 0}\frac{\omega_x (\Pi^{-1}(T^+_x (\Pi
(\zeta),r)))}{r^{1+\xi}}=
    \liminf_{r\searrow 0}\frac{\omega_x (\Pi^{-1}(T^-_x
    (\Pi (\zeta),r)))}{r^{1+\xi}}=+\infty.\]
It is easy to verify if \(K^{(2)}=K\cap (B^{(2)}\times S^1\times
I)\), then \(f\colon  K^{(2)}\to \R\cup\{+\infty\}\) given by
\[f(x,\theta,t)=\min\left\{\liminf_{r\searrow 0}\frac{\nu_x (T^+_x (x+t\hat{\theta},r))}{r^{1+\xi}},\,
    \liminf_{r\searrow 0}\frac{\nu_x (T^-_x
    (x+t\hat{\theta},r))}{r^{1+\xi}}\right\}\]
is a Borel function and so
\[K^{(2)}_\infty=\{\zeta\in  K^{(2)}: f(\zeta)=+\infty\}\]
is a Borel set with \(\omega_x (K^{(2)}\setminus
K^{(2)}_\infty)=0\) for all \(x\in B^{(2)}\). Hence,
\(m^*(K^{(2)}\setminus K^{(2)}_\infty)=0\). Now
\[K_\infty ^{(2)}=\bigcap_{m\in\N}\bigcup_{n\in\N} K^{(2)}_{m,n}\]
where
\begin{align*}
    K^{(2)}_{m,n}=&\left\{ (x,\theta,t)\in K^{(2)}:
    \mbox{if }r\in (0,\tfrac{1}{n}]\mbox{, then}\right.\\
    &\qquad\qquad \left.\min\left\{{\nu_x (T^+
    (x+t\hat{\theta},r))},{\nu_x(T^-
    (x+t\hat{\theta},r))}\right\}>m{r^{1+\xi}}\right\}.
\end{align*}
Thus we can find \(m,n\in\N\) such that \(m^* (K^{(2)}_{m,n})>0\)
and so we can choose a compact set \(K^{(3)}\subseteq
K^{(2)}_{m,n}\) with \( m^* (K^{(3)})>0\). It follows that we can
find a compact set \(B^{(3)}\subseteq \pi_B (K^{(3)})\subseteq
B^{(2)}\) and \(p>0\) such that \(\mu^{(2)} (B^{(3)}) >0\) and for
all \(x\in B^{(3)}\), we have \(\omega_x (K^{(3)})>p\). For \(x\in
B^{(3)}\), let
\[F_x=\Pi (K^{(3)}\cap(\{x\}\times S^1\times I))\subseteq \Gamma_x\cap \Pi^{-1}(U_i)\]
 and notice \(F_x\) is a compact set with \(\nu_x (F_x)=\omega_x (K^{(3)})>p\).
Thus, summarising, we have
\begin{itemize}
    \item \(x\mapsto\nu_x\) is a Borel measurable function on
    \(B^{(3)}\),
    \item  \(B^{(3)}\subseteq \pi_B (K^{(3)})\) is compact with \(\mu^{(2)} (B^{(4)})>0\),
    \item for \(x\in B^{(4)}\), \(\nu_x (F_x)>p\),
    \item for \(x\in B^{(4)}\), \(0<r\leq 1/n\) and
    \(u\in F_x\subseteq \Gamma_x\),
    \[\min\left\{{\nu_x (T^+
    (u,r))},{\nu_x(T^-
    (u,r))}\right\}>m{r^{1+\xi}}.\]
\end{itemize}

Thus we have found a compact set \(B^{(3)}\subseteq E_i\subseteq
B^{(0)}\), a compact set \(\bar{U}_i\cap\Gamma\subseteq \Gamma\),
a non-zero Radon measure \(\mu\) and  a constant \(c>0\) such
that:
\begin{enumerate}
    \item for \(x\in B^{(3)}\) and \(u\in \bar{U}_i\cap\Gamma\),
    \( |u-x|\geq \tfrac{99}{100}d(B,\Gamma)\geq 99\diam (B)\);
    \item \(\mu^{(2)} (B^{(3)}) >0\);
    \item for \(x\in B^{(3)}\) and \(0<r\leq 1\),
    \[\mu B(x,r)\leq r^t;\]
    \item for \(x\in B^{(3)}\), there is a Radon measure \(\nu_x\) and a
    compact set \(F_x\subseteq {U}_i\cap\Gamma_x\cap B_{i} \) such that
    \begin{enumerate}
%        \item \(E_x\subseteq \Gamma_x\cap B_{ij}\);
        \item \(\nu_x (F_x)>p\);
        \item for \(0<r\leq 1\) and \(u\in \Gamma\),
        \[\nu_x B(u,r)\leq r^s ;\]
        \item for \(0<r\leq n^{-1}\) and \(u\in F_x\),
        \[\min\{\nu_{x}(T^+_x(u,r)),\, \nu_{x}(T^-_x(u,r))\}>mr^{1+\xi}.\]
    \end{enumerate}
\end{enumerate}

If \(x,y\in B^{(3)}\), then \(|x-y|\leq \tfrac{1}{100}\dist
(B,\Gamma )\). So let \(d_- =(1/50)\dist (B^{(0)},\Gamma ) \) and
\(d_+=\diam B^{(0)}+\diam\Gamma + \dist (B^{(0)},\Gamma ) \). By
rescaling if necessary, we can assume that \(d_+\leq 1\).

Let \(A, B\subseteq \R^2\) be compact and suppose that \(\psi\in
(0,1/2)\) is such that for \(u\in (A\cap F_x)\cup (B\cap  F_y) \),
\[|\ip{(u-x)^\wedge}{(u-y)^\wedge}|\in [1/ {2},1-\psi ].\]
Then all the hypotheses of Propositions~\ref{prop-geom1}
and~\ref{prop-mass} are satisfied (after suitable relabelling) and
so, for \(u\in A\cap F_x\), \(v\in B\cap F_y\) and \(0<\rho\leq
\dOne\psi^{\frac{1}{2}\frac{1}{s-1-\xi}}\), we find
 that
\[\nu_y ( A\cap F_y\cap B(u,\rho))\leq \dThree \psi^{-\frac{1}{2}\frac{s-1}{2+\xi-s}}
    \rho^{\frac{1+\xi}{2+\xi-s}}\]
and for \(0<\rho\leq \dFour\psi^{\frac{1}{2}\frac{1}{s-1-\xi}}\)
\begin{align}
    \lefteqn{(\nu_x\otimes\nu_y) \left( (F_x\times F_y)\cap(A\times
    B)\cap\{
    (u,v): |u-v|\leq \rho\}\right)}\qquad \qquad\qquad \qquad \label{eqn-mass}\\
    \nonumber
    &\leq \dEight \arcdiam_{\frac{1}{2}(x+y)} (A\cap F_x\cap B(F_y\cap B,\rho))
    (\psi^{-\frac{1}{2}}\rho)^{\frac{s+\xi}{2+\xi-s}}.
\end{align}

\subsection{Energy estimate}

We now pull all our estimates together and explicitly calculate
the $d$-energy of the measure \(\nu\) given by
\[\nu (E)=\int\!\nu_x|_{F_x} (E)\, d\mu|_{B^{(3)}} (x) \mbox{ for  Borel
    \(E\subseteq \R^2\)}\]
and
\[\nu (A)=\inf\{\nu (E): A\subseteq E\mbox{ and }E\mbox{ is Borel}\}
\mbox{, for non-Borel \(A\)}.\]

On noting that for Borel sets \(E\),
\[\nu(E)=\int\omega_x (\Pi^{-1}(E)\cap K^{(3)})\, d\mu|_{B^{(3)}}(x),\]
it is straightforward to verify that \(\nu\) is a Radon measure.
Note that for \(\tau>0\)
\[\int |u-v|^{-\tau} d(\nu \times \nu)(u ,v)
    =\int_{B^{(3)}\times B^{(3)}}\int_{F_x\times F_y} |u-v|^{-\tau}
    d(\nu_x\times\nu_y) (u,v)d(\mu\times\mu)(x,y).\]

Hence, as our choice of \(d\) guarantees  that \(I_{d}=\int |u
-v|^{-d} d(\nu \times \nu)(u ,v)=+\infty\),
\begin{equation}
    \label{eqn-energy}
    \int_{B^{(3)}\times B^{(3)}}\int_{F_{x} \times F_{y}}
    |u-v|^{-d}\,
    d(\nu_{x}  \otimes \nu_{y}) (u,v) {d} (\mu \otimes
    \mu) (x,y)=+\infty.
\end{equation}

Fix \(x\not=y \in B^{(3)}\). In order to reduce writing, we
translate so that \(\tfrac{1}{2}(x+y)=0\) and let \(a=y\), so
\(|x-y|=2|a|\).

Using Fubini's theorem, we find
\begin{align*}
    \lefteqn{\int_{F_{x}\times F_{y}}|u-v|^{-d}\,
    d(\nu_{x}\otimes\nu_{y})(u,v)} & \\
    &= \int_0^\infty (\nu_{x}\otimes\nu_{y})\left( \{(u,v)\in
    F_x\times F_y: |u-v|^{-d}\geq r\}\right)\, dr \\
    &= d\int_0^\infty \rho^{-d-1}(\nu_{x}\otimes\nu_{y})\left( \{(u,v)\in
    F_x\times F_y: |u-v|\leq \rho\}\right)\, d\rho \\
    &= d\int_{F_x}\int_0^\infty \rho^{-d-1}\nu_y (F_y\cap
    B(u,\rho))\, d\rho\, d\nu_x (u) .
\end{align*}

Let
\[A_0^+ =\{w\in A(0,d_-,d_+):\ip{w}{a^\perp}\geq |\ip{w}{a}|\},\]
\[A_0^- =\{w\in A(0,d_-,d_+):\ip{w}{a^\perp}\leq -|\ip{w}{a}|\},\]
and for \(m,n\in\{0,1\}\) and \(i\in\N\), set
\begin{align*}
    A_{i}^{mn} =&\{w\in A(0,d_-,d_+): |\ip{w}{a^\perp}/\ip{w}{a}|\in
        [2^{-i},2^{1-i}],\, \\
        &\quad\qquad\qquad (-1)^n\ip{w}{a^\perp }>0\mbox{ and
        }(-1)^m\ip{w}{a}>0\},
\end{align*}
noting that
\[\arcdiam_0 (A_i^{mn})\preceq 2^{-i}.\]

Observe that if \(w\) is in \(A_i^{mn}\),  then, by
lemma~\ref{lem-angle},
\[\frac{1}{2}\leq \frac{\ip{w-a}{w+a}}{|w-a||w+a|}\leq 1-
\frac{9}{17d_+^2}(|a|2^{-i})^2,\] and if  \(w\in A_0^+\cup
A_0^-\), then
\[\frac{1}{2}\leq \frac{\ip{w-a}{w+a}}{|w-a||w+a|}\leq 1- \frac{9}{17d_+^2}|a|^2.\]

 For \(i\in\N\cup\{0\}\), set
\(\psi_i =\left(\tfrac{2}{5}\tfrac{d_-}{d_2}|a|2^{-i}\right)^2\)
and let \(\rho_i= \dFour \psi_i^{\frac{1}{2}\frac{1}{s-1-\xi}}\).
Observe that, since \(\rho_i\leq \frac{2}{5}d_- 2^{-i}\),
 if
\(u\in A_i^{mn}\) and \(|u-v|\leq \rho_i\) with \(v\in
A(0,d_-,d_+)\), then \(v\in A_{i-1}^{mn}\cup A_i^{mn}\cup
A_{i+1}^{mn}\). Similarly, if \(u\in A_0^+\cup A_0^-\), and \(v\in
A(0,d_-,d_+)\) with \(|u-v|\leq \rho_0\), then \(v\in A_1^{mn}\)
for some choice of \(m\) and \(n\).

\begin{figure}
\begin{center}
\includegraphics[width=0.7\textwidth, clip=true, angle=90, keepaspectratio=true]{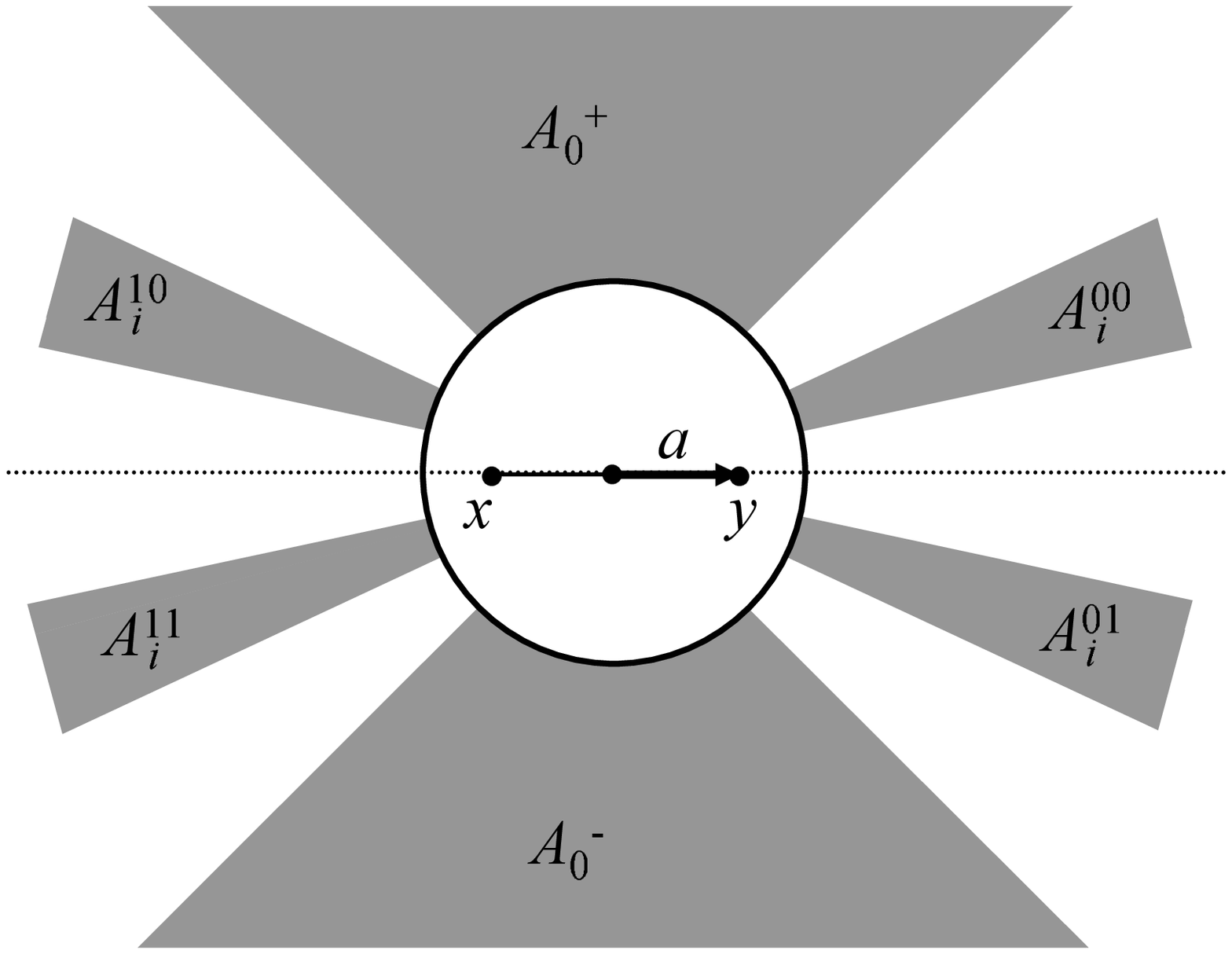}
\end{center}
    \caption{\(A_i^{mn}\).}\label{fig-final}
\end{figure}

Writing \(f(\rho)=\rho^{-d-1}\nu_y   (F_y\cap B(u,\rho))\), we let
\(I_0^{\pm}= \int_{F_x\cap A_0^\pm}\int_0^\infty f(\rho )\, d\rho
d\nu_{x} (u)\) and \(I_i^{mn}=\int_{F_x\cap A_i^{mn}}\int_0^\infty
f(\rho)d\rho\,d\nu_x (u)\).

We must estimate
\begin{align*}
    \lefteqn{\int_{F_x} \int_0^\infty \rho^{-d-1}\nu_y
        (F_y\cap B(u,\rho))\, d\rho\,d\nu_x (u)  }&\\
        & =\left(\int_{F_x\cap A_0^+}+\int_{F_x\cap A_0^-}+
        \sum_{m,n=0}^1\sum_{i=1}^\infty
        \int_{F_x\cap A_i^{mn}}\right)\int_0^\infty f(\rho)\, d\rho\,d\nu_x (u)\\
        &= I_0^+ +I_0^- +\sum_{m,n=0}^1\sum_{i=1}^\infty I_i^{mn}.
\end{align*}
We can write
\begin{align*}
    I_i^{mn}&=\int_{F_x\cap A_i^{mn}}\left(\int_0^{\rho_i}+\int_{\rho_i}^\infty\right)f(\rho)\, d\rho
    \,d\nu_x (u)\\
    &= \int_{F_x\cap A_i^{mn}\cap B(F_y\cap B(A_i^{mn},\rho_i),\rho_i)}\int_0^{\rho_i}f(\rho)\, d\rho\, d\nu_x (u) +
    \int_{F_x\cap A_i^{mn}}\int_{\rho_i}^\infty f(\rho)\, d\rho\, d\nu_x
    (u)\\
    &= I_{i,1}^{mn}+I_{i,2}^{mn}\mbox{, say.}
\end{align*}

\begin{lem}\label{lem-energy2}
    Suppose that \(V\subseteq A(0,d_-, d_+)\) and \(0<r <1\). Then
    \[ \int_{F_x\cap V}\int_r^\infty f(\rho) \, d\rho \leq cr^{s-d}
    \arcdiam_0 (F_x\cap V)^{s-1},\]
    where \(c\) is a positive constant that depends
    only on \(d_-\), \(d_+\), \(s\) and \(d\).
\end{lem}

In the proof of the lemma, and subsequently, we let \(\preceq\)
denote inequality up to a finite  constant independent of $x$ and
$y$.

\begin{proof}
    Using the crude estimate that for \(u\in F_x\cap V\),
    \( \nu_y (F_y\cap B(u,\rho))\leq
    \min\{1,\,2^s\rho^s\}\) together with Lemma~\ref{lem-massnuxfx},
    we find
    \begin{align*}
    \lefteqn{\int_{F_x\cap V}\int_r^\infty f(\rho) \, d\rho }\\
    &\leq 2^s\left(\frac{1}{d-s}r^{s-d}+
        \frac{1}{d}\right)\nu_x (F_x\cap V)\\
        &\preceq r^{s-d}\arcdiam_0 (F_x\cap V)^{s-1}.
    \end{align*}
\end{proof}

In particular, Lemma~\ref{lem-energy2} implies that
\[I_{i,2}^{mn}\preceq \rho_i^{s-d}\arcdiam_0 (F_x\cap A_i^{mn})^{s-1}
\preceq
|a|^{\frac{s-d}{s-1-\xi}}2^{-\left(\frac{s-d}{s-1-\xi}+s-1\right)i}.\]

In order to estimate \(I_{i,1}^{mn}\),  we use
equation~(\ref{eqn-mass}) (of section~\ref{ssec-decomp}), Fubini's
theorem  and the fact that if \(u\in F_x\cap A_i^{mn}\) and \(v\in
B(u,\rho_i)\cap F_y\), then \(v\in  A_{i-1}^{mn}\cup
A_{i}^{mn}\cup A_{i+1}^{mn}\), to calculate  that, provided
\(\frac{s+\xi}{2+\xi-s}-d >0\),
\begin{align*}
    \lefteqn{I_{i,1}^{mn}}\\
    &=\int_0^{\rho_i}\rho^{-d-1}
    (\nu_x\otimes \nu_y)\{ (u,v)\in (F_x\cap A_{i}^{mn})\times (F_y\cap B(A_i^{mn},\rho))
    :|u-v|\leq \rho\}\, d\rho\\
    &\leq  \dEight \arcdiam_{0} (A_i^{mn}\cap F_x)\int_0^{\rho_i} \rho^{-d-1}
    (\psi_{i+1}^{-\frac{1}{2}}\rho)^{\frac{s+\xi}{2+\xi-s}}\, d\rho\\
    &\preceq \arcdiam_{0} (A_i^{mn}\cap F_x)\psi_{i+1}^{-\frac{1}{2}\frac{s+\xi}{2+\xi-s}}
    \rho_i^{\frac{s+\xi}{2+\xi-s}-d}\\
    &\preceq \arcdiam_{0} (A_i^{mn}\cap F_x)\psi_{i}^{\frac{1}{2}\left(
    \frac{s+\xi-d}{s-1-\xi}\right)}\\
    &\preceq
    |a|^{\frac{s+\xi-d}{s-1-\xi}}2^{-i\left(\frac{2s-1-d}{s-1-\xi}\right)}.
\end{align*}

Combining these estimates for \(I_{i,1}^{mn}\) and
\(I_{i,2}^{mn}\), we deduce that, provided
\(\frac{s+\xi}{2+\xi-s}-d >0\),
\begin{align*}
    I_{i}^{mn} &\preceq |a|^{\frac{s+\xi-d}{s-1-\xi}}
        2^{-i\left(\frac{2s-1-d}{s-1-\xi}\right)}+
        |a|^{\frac{s-d}{s-1-\xi}}(2^{-i})^{\frac{s-d}{s-1-\xi}
    +s-1}\\
    &= |a|^{\frac{s-d}{s-1-\xi}}\left((2^{-i})^{\frac{2s-1-d}{s-1-\xi}}
    +(2^{-i})^{\frac{s-d}{s-1-\xi}
    +s-1}\right).
\end{align*}

Hence
\[\sum_{m,n=0}^1\sum_{i=1}^\infty I_{i}^{mn}\preceq |a|^{\frac{s-d}{s-1-\xi}},\]
provided that
\(\min\left\{\frac{s+\xi}{2+\xi-s}-d,\,\frac{2s-1-d}{s-1-\xi},\,\frac{s-d}{s-1-\xi}
    +s-1\right\}>0 \).
Estimating \(I_0^+\) and \(I_0^-\) in a similar way, we find
\[I_0^++I_0^-+\sum_{m,n=0}^1\sum_{i=1}^\infty I_{i}^{mn}\preceq
    |a|^{\frac{s-d}{s-1-\xi}},\]
provided that
\(\min\left\{\frac{s+\xi}{2+\xi-s}-d,\,\frac{2s-1-d}{s-1-\xi},\,\frac{s-d}{s-1-\xi}
    +s-1\right\}>0\).
Hence, provided that
\(\min\left\{\frac{s+\xi}{2+\xi-s}-d,\,\frac{2s-1-d}{s-1-\xi},\,\frac{s-d}{s-1-\xi}
    +s-1\right\}>0\),
\[\int_{F_{x}\times F_{y}}|u-v|^{-d}\,
    d(\nu_{x}\otimes\nu_{y})(u,v)
    \preceq |x-y|^{-\frac{d-s}{s-1-\xi}}.\]
Thus, if
\(\min\left\{\frac{s+\xi}{2+\xi-s}-d,\,\frac{2s-1-d}{s-1-\xi},\,\frac{s-d}{s-1-\xi}
    +s-1\right\}>0\), then
\[+\infty=I_d(\nu)\preceq I_{\frac{d-s}{s-1-\xi}}(\mu)\]
and this gives a contradiction if \(\frac{d-s}{s-1-\xi}<t\), the
dimension of \(\mu\). Since \(0<\xi<s-1\) is arbitrary, it follows
that if
\(s>\max\left\{\frac{1}{2}(d+1),\frac{2d}{d+1},\frac{1}{2}+\sqrt{d-\frac{3}{4}}\right\}=
\frac{1}{2}+\sqrt{d-\frac{3}{4}}\), then \(t\leq \frac{d-s}{s-1}\)
and Theorem~\ref{thmresult} (and hence Theorem~\ref{thmresult2})
follows.

\end{document}